\documentclass{compositio}
\usepackage[all]{xy}

\usepackage{graphicx}
\usepackage{latexsym,amsthm,amssymb,amscd,youngtab,stmaryrd}
\usepackage{enumerate}
\usepackage{multicol}
\usepackage{amscd}
\usepackage{amsmath}

\usepackage{floatflt}
\numberwithin{equation}{subsection}
\numberwithin{figure}{section}

\newcommand{\test}{\it}
\newcommand{\Yv}{\Yvcentermath1}

\newtheorem{theorem}{Theorem}[subsection]
\newtheorem{lemma}[theorem]{Lemma}
\newtheorem{remark}[theorem]{Remark}
\newtheorem{prop}[theorem]{Proposition}
\newtheorem{cor}[theorem]{Corollary}
\newtheorem{ex}[theorem]{Example}
\newtheorem{conjecture}[theorem]{Conjecture}
\newtheorem{definition}[theorem]{Definition}

\newtheorem{THEOREM}{Theorem}
\newtheorem{corr}{Corollary}

\newcommand{\cA}{{\mathcal A}}
\newcommand{\cB}{{\mathcal B}}
\newcommand{\cC}{{\mathcal C}}
\newcommand{\cD}{{\mathcal D}}
\newcommand{\cE}{{\mathcal E}}
\newcommand{\cF}{{\mathcal F}}
\newcommand{\cG}{{\mathcal G}}
\newcommand{\cH}{{\mathcal H}}
\newcommand{\cI}{{\mathcal I}}
\newcommand{\cK}{{\mathcal K}}
\newcommand{\cL}{{\mathcal L}}

\newcommand{\cO}{{\mathcal O}}
\newcommand{\cP}{{\mathcal P}}

\newcommand{\cQ}{{\mathcal Q}}
\newcommand{\cS}{{\mathcal S}}

\newcommand{\cU}{{\mathcal U}}

\newcommand{\cZ}{{\mathcal Z}}

\newcommand{\mg}{\mathfrak{g}}
\newcommand{\mh}{\mathfrak{h}}
\newcommand{\mb}{\mathfrak{b}}

\newcommand{\mm}{\mathfrak{m}}
\newcommand{\mn}{\mathfrak{n}}

\newcommand{\mR}{\mathbb{R}}
\newcommand{\mV}{\mathbb{V}}

\newcommand{\mC}{\mathbb{C}}
\newcommand{\mH}{\mathbb{H}}

\newcommand{\mZ}{\mathbb{Z}}

\newcommand{\la}{\lambda}
\newcommand{\oneone}{\stackrel{1:1}\leftrightarrow}

\newcommand{\HOM}{\operatorname{Hom}}

\newcommand{\EXT}{\operatorname{Ext}}

\newcommand{\END}{\operatorname{End}}

\newcommand{\p}{\mathfrak{p}}
\newcommand{\MOD}{\operatorname{mod}}

\newcommand{\op}{\operatorname}

\newcommand{\KER}{\operatorname{ker}}

\newcommand{\COKER}{\operatorname{coker}}

\newcommand{\DIM}{\operatorname{dim}}

\newcommand{\inj}{\hookrightarrow}

\newcommand{\surj}{\mbox{$\rightarrow\!\!\!\!\!\rightarrow$}}
\newcommand{\iso}{\tilde{\rightarrow}}

\newcommand{\CAN}{\operatorname{can}}
\newcommand{\ID}{\test{id}}

\begin{document}
\title[Parabolic category $\cO$, Springer fibres and Khovanov homology]{Parabolic category $\cO$, Perverse sheaves on Grassmannians, Springer fibres and Khovanov homology}

\author{Catharina Stroppel}
\email{c.stroppel\symbol{64}maths.gla.ac.uk}
\address{Department of Mathematics, University of Glasgow (UK).}
\thanks{This work was supported by EPSRC grant 32199}
\subjclass{16S99, 17B10, 14M15, 57M27, 20C30, 20G05,14M17}

\begin{abstract}
For a fixed parabolic subalgebra $\p$ of $\mathfrak{gl}(n,\mC)$ we prove that the centre of the principal block $\cO_0^\p$ of the parabolic category $\cO$ is naturally isomorphic to the cohomology ring $H^*(\cB_\p)$ of the corresponding Springer fibre. We give a diagrammatic description of $\cO_0^\p$ for maximal parabolic $\p$ and give an explicit isomorphism to Braden's description of the category $\op{Perv}_B(G(k,n))$ of Schubert-constructible perverse sheaves on Grassmannians. As a consequence Khovanov's algebra $\cH^n$ is realised as the endomorphism ring of some object from $\op{Perv}_B(G(n,n))$ which corresponds under localisation and the Riemann-Hilbert correspondence to a full projective-injective module in the corresponding category $\cO_0^\p$. From there one can deduce that Khovanov's tangle invariants are obtained from the more general functorial invariants in \cite{StDuke} by restriction.
\end{abstract}
\maketitle

\tableofcontents

\section*{Introduction}

Let $G=\op{GL}(n,\mC)$ be the general linear group with subgroup $B$ given by all invertible upper triangular matrices.
Let $\mg=\mathfrak{gl}(n,\mC)$ and $\mb$ be their Lie algebras and let $W$ be the Weyl group, so $W=S_n$.
Let $\mu=(\mu_1, \mu_2,\ldots ,\mu_r)$ be positive integers summing up to $n$.
Then we have the parabolic subalgebra $\p\supseteq \mb$ of $\mg$ with Levi subalgebra
$\mathfrak{gl}_{\mu_{1}}\oplus \mathfrak{gl}_{\mu_{2}}\oplus\ldots \oplus \mathfrak{gl}_{\mu_{r}}$,
Weyl group $W_\p=S_{\mu_1}\times  S_{\mu_2}\times\ldots\times S_{\mu_r}$, and $P$ the corresponding parabolic subgroup of $G$.
Let $x_\mu\in\mg$ be a nilpotent element whose Jordan Normal Form has blocks of size $\mu_i$, $1\leq i\leq r$. Let $u_\mu=\operatorname{Id}+x_\mu$ be the corresponding unipotent element. Let  $\cB=G/B$, the variety of full flags in $\mC^n$. Associated to $\mu$ we have the partial flag variety $G/P$. We also have the variety $\cF^{u_\mu}$ of $u_\mu$-fixed points on the other hand: $\cF^{u_\mu}$ is the Springer fibre associated with $\p$; we denote it by $\cB_\p$. \\

Let $\cO_0^\p$ be the principal block in the category of
highest weight modules for $\mg$ which are locally finite
for $\p$. If $\p=\mb$ then $\cO_0^\p$ is the principal
block of the ordinary Bernstein-Gelfand-Gelfand category
$\cO$. The category $\cO_0^\p$ is equivalent to
$\MOD-A^\p$, the category of finitely generated modules
over the (finite dimensional) endomorphism algebra $A^\p$
of a minimal projective generator of $\cO_0^\p$; it is also
equivalent to the category of perverse
sheaves on $G/P$, constructible with respect to the Schubert stratification (via localisation and the Riemann-Hilbert
correspondence).

\subsubsection*{Centres and Springer fibres}
The first result of the paper generalises Soergel's results from \cite{Sperv}, confirms \cite[Conjecture 3]{KhoSpringer} and gives an explicit description of the centre $Z(A^\p)$ of $A^\p$:

\begin{THEOREM}
\label{theoone}
There is a canonical isomorphism of algebras $H^*(\cB_\p)\cong\cZ(A^\p)$.
\end{THEOREM}
In particular, up to isomorphism, $\cZ(A^\p)$ only depends on the parts of $\mu$, not on the order in which they appear.

This theorem was independently proved by Jonathan Brundan (\cite{Brundan2})
using different techniques. His approach also works for
singular blocks and provides an explicit description of these centres as
quotients of polynomial rings.

The cohomology of Springer fibres  $H^*(\cB_\p)$ was used by Springer (\cite{Springer}) to construct the irreducible representations of the symmetric group $S_n$. In particular, he defined an $S_n$-action on $H^*(\cB_\p)$. From the isomorphism above we get an induced $S_n$-action on $\cZ(A^\p)$. In Section~\ref{centreacts} we give a functorial interpretation of this $S_n$-action on  $\cZ(A^\p)$ as follows. Let $B_W$ be the underlying braid group. Jantzen's translation functors can be used to define a (weak) braid group action on the bounded derived category of $\cO^\p_0$ (see e.g. \cite{StDuke}, \cite{BFK}, \cite{Rouquier}).
The resulting functors are the derived functors of Irving's shuffling functors (see Section \ref{shuffle}).
Since these functors are tilting functors, they induce a braid group action on the centre of $\cO^\p_0$, hence on $\cZ(A^\p)$. Now the natural map $can:\cZ(A^\mb)\rightarrow\cZ(A^\p)$ is $B_W$-equivariant (Theorem~\ref{centreb}) and if $\mb=\p$ then the braid group action factors through an action of $S_n$ (Lemma~\ref{action}). A very recent result of Brundan (\cite{Brundan}) says that $can$ is surjective, hence there is an $S_n$-action on $\cZ(A^\p)$ as well. This is the $S_n$-action we are looking for. Together with the theorem above and the main result from either \cite{CoPr} or \cite{Tanisaki} it follows that $\cZ(\cO^\p_0)\cong\mC[W]\otimes_{\mC[W_\p]}\mC_{triv}$ as $W$-module. It also shows that the dimension of the centre stays invariant under deformations of $A^\p$. \\

The main idea of the proof of Theorem~\ref{theoone} is as
follows: From Soergel's Endomorphismensatz and Struktursatz
\cite{Sperv} we get an isomorphism
$H^*(\cB)\cong\cZ(A^\mb)$ of rings. On the other hand we
have the restriction map $\cZ(A^\mb)\rightarrow\cZ(A^\p)$.
We first show that the kernel of the canonical map
$H^*(\cB)\rightarrow H^*(\cB_\p)$ is contained in the
kernel of $\cZ(A^\mb)\rightarrow\cZ(A^\p)$ using
deformation theory (following \cite{Sperv}). This is based
on the results of \cite{Tanisaki} and a handy description
of  $H^*(\cB_\p)$ as a quotient of $S(\mh)$ along the lines
of \cite{GP}. To show that the induced map $\Phi_\p:
H^*(\cB_\p)\rightarrow \cZ(A^\p)$ is {\it injective} it is
enough to show that it is injective on its socle
(considered as an $H^*(\cB)$-module). The main idea here is
that the top degrees of $H^*(\cB_\p)$ and $\cZ(A^\p)$
coincide (Lemma~\ref{top}).  This will be used to show that
$\Phi_\p$ is non-zero when restricted to the socle
(Proposition~\ref{Irving}, Proposition~\ref{ll}), and even
$S_n$-equivariant onto its image. Then we use the fact that
the socle of $H^*(\cB_\p)$ is an irreducible $S_n$-module
and get the injectivity. In Theorem~\ref{SpringerIrr} we
show that the induced injective map defines an isomorphism
of $S_n$-modules on the top degree parts
\begin{eqnarray*}
  H^*(\cB_\p)_{top}\iso\cZ(A^\p)_{top}.
\end{eqnarray*}
This follows on the one hand from Springer's construction of irreducible
$S_n$-modules, and on the other hand from the categorification of irreducible
$S_n$-modules using projective-injective modules in $\cO_0^\p$ obtained in
\cite{KMS}. In the maximal parabolic case we give an alternative proof for the
injectivity by a deformation argument, since the algebra $A^\p$ can be replaced
be a symmetric subalgebra with the same centre (Section~\ref{max}) and the
deformation ring is a principal ideal domain. For the applications we have in
mind (see below) the maximal parabolic case is enough. As far as we see
deformation methods are not sufficient to prove the surjectivity in
Theorem~\ref{theoone}. Instead, some `external' information is needed which is
obtained in \cite{Brundan} from the representation theory of cyclotomic Hecke
algebras.

\subsubsection*{Connection to Khovanov homology}
Theorem~\ref{theoone} together with \cite{CoPr} and \cite{Tanisaki} provide an explicit description of $\cZ(A^\p)$ - so, we would like to have an explicit description of $A^\p$ as well. In general, this seems to be ambitious, but  in case $\p$ is a maximal parabolic subalgebra it has been achieved by Braden (\cite{Braden}) using perverse sheaves on Grassmannians. However, the description there is difficult to use for explicit calculations. Moreover, the Koszul grading of $A^\p$ (defined in \cite{BGS}) is not visible. Therefore, we consider the situation of \cite{Braden} again and first remark that any indecomposable projective $A^\p$-module has a commutative endomorphism ring (Proposition~\ref{comm}).
Later on we deduce that each of these endomorphism rings is of the form $\mC[X]/(X^2)^{\otimes k}$ for some $k\in\mZ_{\geq0}$. In Corollary~\ref{corgrad} we explain how $A^\p$ becomes a graded algebra using the description of \cite{Braden}. The intriguing result is however Theorem~\ref{graphical} which gives a purely graphical description of Braden's algebra $A_{m,m}$ very similar to Khovanov's approach (see e.g. \cite{KhoSpringer}) which we will describe below. \\

The crucial fact behind Theorem~\ref{theoone} and its proof is the existence of
a bijection between the isomorphism classes $\op{PrInj}(\p)$ of indecomposable
projective-injective modules in $\cO_0^\p$ and the irreducible components in
$\cB_\p$. Let us consider the case where $n=2m$ for some $m\in\mZ_{>0}$ and
$\mu=(m,m)$. In this case the irreducible components of $B_\p$, and hence the
isomorphism classes $\op{PrInj}(\p)$, are in bijection to $I$, the set of
crossingless matchings of $2m$ points.
Let $\{T(i)_{2m}\}_{i\in I}$ be a complete minimal set of representatives of $\op{PrInj}(\p)$ and $T_{2m}:=\oplus_{i\in I} T(i)_{2m}$. In \cite{Khotangles}, a finite dimensional $\mC$-algebra $\cH^m$ was introduced whose primitive idempotents are naturally indexed by crossingless matchings of $2m$ points. These algebras were used to define the famous Khovanov homology which gives rise to an invariant of tangles and links. It is known (\cite{KhoSpringer}) that the centre $\cZ(\cH^m)$ of $\cH^m$ is isomorphic to $H^*(\cB_\p)$.\\

In Theorem~\ref{arbeit} and Proposition~\ref{surjectivity} we verify \cite[Conjecture 2.9 (a)]{StTQFT} which is a stronger version of the conjectures formulated in \cite{Braden} and \cite{KhoSpringer}:

\begin{THEOREM}
\label{conjTQFT}
For any natural number $m$ there is an isomorphism of algebras
\begin{eqnarray*}
\END_\mg(T_{2m})\iso\cH^m.
\end{eqnarray*}
\end{THEOREM}

Hence Khovanov's algebra $\cH^m$ is a subalgebra of $A^\p$, where $\p$ is the parabolic subalgebra of $\mathfrak{gl}_{2m}$ corresponding to the decomposition $2m=m+m$.

\begin{corr}
  There is an isomorphism of rings $\cZ(\END_\mg(T_{2m}))\iso\cZ(\cH^m)$
\end{corr}

With \cite[Theorem 3]{KhoSpringer} we therefore have an alternative proof of Theorem~\ref{theoone} in this special situation (purely based on \cite{Braden}) which implies \cite[Conjecture 2]{KhoSpringer}. \\

As an application of Theorem~\ref{conjTQFT} one can deduce that Khovanov's tangle invariants are nothing else than restrictions of the functorial invariants from \cite{StDuke} (see \cite[Conjecture 2.9(b)]{StTQFT} for a precise statement). Since the proof is lengthy, this part will be presented in a subsequent paper.

\subsubsection*{Diagrammatic description of Braden's algebra}
Theorem~\ref{conjTQFT} will be a direct consequence of our diagrammatic
description of Braden's algebra $A_{m,m}$ in the case $\mu=(m,m)$. The
primitive idempotents of $A_{m,m}$, or equivalently the isomorphism classes of
indecomposable projective modules of $\cO_0^\p$, are in bijection to the
shortest coset representatives of $S_m\times S_m \backslash S_{2m}$: by
permuting the entries, the symmetric group $S_{2m}$ acts transitively on the
set of $\{+,-\}$-sequences of length $2m$ with exactly $m$ pluses and $m$
minuses. Since the sequence $\sigma_{dom}=(+,\ldots, +,-, \ldots, -)$ has
stabiliser $S_m\times S_m$ we get a bijection between the primitive idempotents
of $A_{m,m}$ and the set $\cS(m)$ of $\{+,-\}$-sequences of length $2m$ with
exactly $m$ pluses and $m$ minuses (Proposition~\ref{bijall}). The isomorphism
class of the projective generalised Verma module in $\cO^\p_0$ is mapped to
$\sigma_{dom}$ under this bijection. For $m=1$ we have the sequence $(+,-)$
corresponding to the projective Verma module and the sequence $(-,+)$
corresponding to the `antidominant projective module'.

We want to associate a cup-diagram to each isomorphism class of indecomposable
projective modules. To do so we have to make the $\{+,-\}$-sequences longer.
Putting $m$ minuses in front of a sequence from $\cS(m)$ and $m$ pluses
afterwards we obtain a distinguished set of $\{+,-\}$-sequences of length $4m$
with exactly $2m$ pluses. Connecting successively each minus with an orphaned
neighboured plus to the right we obtain a collection of crossingless matchings
of $4m$-points.  In this way we associate to each primitive idempotent $a$ of
$A_{m,m}$ a cup diagram/crossingless matching of $4m$ points. For the sake of
argument in this introduction we number the points from $1$ to $4m$. In the
case $m=1$ for example, the two sequences $(-,+)$ and $(+,-)$ of length $2$
from above become the sequences $(-,-,+,+)$ and $(-,+,-,+)$ of length $4$ and
we associate the crossingless matchings depicted in Section~\ref{TL},
Figure~\ref{fig:cups}.

To a pair $(a,b)$ of two primitive idempotents we obtain a collection of circles as in \cite{Khotangles}, namely by putting one crossingless matching upside down on top of the other (see Section~\ref{TL} Figure~\ref{fig:gluing} for $m=1$).

The fundamental difference to \cite{Khotangles} is that we
additionally introduce a colouring of these circles
indicating the position of a circle (Section~\ref{TL},
Figure~\ref{fig:colours}). If a circle connects only points
in the interval $[m+1,3m]$ then the circle is black. If a
circle passes either through at least two points in $[1,m]$
or at least two points in $[3m+1,4m]$ then it is red. In
all other cases it is green.

The principle idea is that we fix for each allowed colour (black, red, green) a 2-dimensional TQFT: red circles correspond to the trivial Frobenius algebra, green circles correspond to the one dimensional Frobenius algebra, and black circles correspond to the Frobenius algebra $\mC[X]/(X^2)$. In Section~\ref{enlarged} we will combine these three TQFTs to define an algebra $\cK^m$. If we restrict to idempotents such that only black circles occur then we are exactly in the situation of \cite{Khotangles} and we obtain $\cH^m$ naturally as a subalgebra of $\cK^m$. However, the colouring carries all the additional information to give a graphical description of Braden's algebra $A_{m,m}$ (see Theorem~\ref{graphical}):

\begin{THEOREM}
\label{Theothree}
  For any $m\in\mZ_{>0}$ there is an isomorphism of algebras $$\cE: A_{m,m}\cong\cK^m.$$
\end{THEOREM}

To determine the dimension of the homomorphism space between two indecomposable
projective modules it is enough to take the corresponding two cup diagrams, one
upside down on top of the other, and count the numbers of circles for each
colour. The dimension of the morphism space is then zero if there is a red
circle, otherwise two to the power of the number of black circles.
Taking the dimension of these homomorphism spaces is a natural extension of the
categorical version of the $S_n$-invariant bilinear form defined on irreducible
$S_n$-modules described in \cite{KMS}.

The homomorphism space between two indecomposable projective modules $P$ and
$Q$ carries a natural $\mZ$-grading induced from the Koszul grading introduced
in \cite{BGS}. It turns out that, up to a shift, the Poincar{\'e} polynomial
agrees with the Intersection Theory Poincar{\'e} polynomial associated to the
intersection of the corresponding two irreducible components of the associated
Springer fibre (Theorem~\ref{PrInjSpringer}).

On the other hand $\cK^m$ carries a natural $\mZ$-grading
induced from the $\mZ$-grading on $\mC[X]/(X^2)$, where $X$
has degree two. The isomorphism from
Theorem~\ref{Theothree} induces a grading of $A_{m.m}$
which we show is the Koszul grading (Corollary
~\ref{corgrad}). It follows in particular, that the arrows
in the Ext-quiver of $A_{m,m}$ are given by Braden's
relation $\leftrightarrow$ (Corollary~\ref{Ext}).

\subsubsection*{Plan} The paper starts by recalling basics from Category $\cO$ and
its deformation theory in Section~\ref{defo}.
Section~\ref{centreacts} contains general facts about braid
group actions on the centres of the categories we are
interested in. Starting from Section~\ref{typeA} we will
only consider Type $A$, the Lie algebra $\mathfrak{gl}_n$.
In Section~\ref{typeA} we explain the connection between
the centres of blocks of category $\cO$ and the cohomology
of the Springer fibres. Section~\ref{TL} contains the
connection with Braden's and Khovanov's work. We tried to
make this part accessible without the Lie theoretic
background from the previous Sections. We abbreviate
$\otimes_\mC$ as $\otimes$ and
$\operatorname{dim}=\operatorname{dim}_\mC$ denotes the
dimension of a complex vector space.

\subsubsection*{Acknowledgement} I am grateful to Tom Braden for all his explanations -
the results appearing in Section \ref{SectionBraden} are built on his ideas,
and to Toshiyuki Tanisaki - crucial steps towards the result of this paper were
obtained during my visit to Osaka in December 2005. Thanks to Henning Haahr
Andersen for very useful remarks on a previous version, and Volodymyr Mazorchuk
and Wolfgang Soergel for useful discussions on the way. I would like to thank
Jonathan Brundan for letting me know about his results which improve
Theorem~\ref{theosingular}. Thanks to Iain Gordon for always sharing his ideas,
for his constant support and many useful comments and suggestions. I would in
particular like to thank the referees for extremely useful comments and
pointing out a gap in the arguments.

\section{Preliminaries}
Let $G$ be a complex reductive simply connected algebraic group with a chosen
Borel subgroup $B$ and maximal torus $T$. Let $\mg$ be the corresponding
reductive complex Lie algebra, with $\mb\supset\mh$ the Lie algebras of $B$ and
$T$ respectively. For any Lie algebra $\mathfrak{l}$ let $\cU(\mathfrak{l})$ be
its universal enveloping algebra. We abbreviate $\cU=\cU(\mg)$ and denote by
$\mathfrak{Z}$ the centre of $\cU$. For any finite dimensional complex vector
space $V$ let $S(V)$ be the algebra of polynomial functions on $V^*$,
especially $S:=S(\mh)=\cU(\mh)$. We denote by $\triangle\subseteq R^+\subset R$
the set of simple roots, positive roots and all roots. Let $W$ be the Weyl
group and $X=X(R)$ the integral weight lattice. The Weyl group acts naturally
on $\mh^*$; this action is denoted by $(w,\la)\mapsto w(\la)$, for $w\in W$,
$\la\in\mh^*$. We also have
 the so-called `dot-action' given by $(w,\la)\mapsto
w\cdot\la:=w(\la+\rho)-\rho$, where $\rho$ is the half-sum of
positive roots. For a root $\alpha\in R$ we denote by
$\check{\alpha}$ the corresponding coroot with the evaluation
pairing $\langle,\rangle$. In this paper, a weight $\la\in\mh^*$
is called {\it
  dominant}, if $\langle\la+\rho,\check{\alpha}\rangle\notin\{-1,-2,\ldots\}$ for any $\alpha\in R^+$. Let $\mh^*_{dom}$ be the set of dominant
weights.

If $\pi\subseteq\triangle$ then there is a corresponding parabolic subalgebra
$\p_\pi=\mg_\pi\oplus\mh^\pi\oplus\mn^\pi$ of $\mg$, where $\mg_\pi$ is
semisimple with simple roots $\pi$ and Cartan $\mh_\pi$, and
$\mh^\pi=\cap_{\alpha\in\pi}\op{ker}{\alpha}$. Denote by $p_\pi:\mh\rightarrow
\mh^\pi$ the projection along $\mh_\pi$ as well as the induced restriction
morphism $p_\pi:S(\mh)\rightarrow S(\mh^\pi)$.

Let $W_{\p_{\pi}}$ be the parabolic subgroup of $W$ associated with
$\pi$. If $\pi=\emptyset$ then $W_{\p_\pi}$ is trivial. We denote by $W^{\p_\pi}$ the set of shortest coset
representatives in $W_{\p_{\pi}}\backslash W$ with respect to the
Coxeter length function $l$. Let $w_0\in W$ be the longest
element, $w_0^{\p_\pi}$ the longest element in $W_{\p_{\pi}}$, and
finally $[w^{\p_\pi}]$ the representative of the longest element of
$W_{\p_{\pi}}\backslash W$.

\section{Deformation}
\label{defo} Fix now some $\pi\subseteq\triangle$. To simplify notation we will
leave out the index $\pi$ most of the time. In particular, $\p=\p_\pi$.
Consider the following commutative algebra
\begin{eqnarray*}
  T=T^\pi:=S(\mh^\pi)_{(0)}=\left\{\frac{f}{g}\mid f,g\in S(\mh^\pi),
  g(0)\not=0\right\},
\end{eqnarray*}
the localisation of $S(\mh^\pi)$ at the augmentation ideal, with maximal ideal $\mm$.

 Let $T'$ be a {\it $\pi$-deformation algebra}, that is an associative unitary noetherian commutative $T$-algebra. The structure morphism $\varphi:T\rightarrow T'$, induces a $\cU(\p)$-module structure on $T'$ via the
composition $\p_\pi\stackrel{p}{\longrightarrow} \mh^\pi\stackrel{i}
  \longrightarrow S(\mh^\pi)\stackrel{i'}\rightarrow T\stackrel{\varphi}{\longrightarrow} T',$
where $p$ is the canonical projection, and $i$, $i'$ the canonical inclusions. We have the set $X_\pi$ of $\pi$-integral weights and the set $X_\pi^+$ of ($\pi$-)admissible weights, defined as follows:
\begin{eqnarray*}
  X_\pi&:=&\{\la\in\mh^*\mid \langle \la,\check{\alpha}\rangle\in\mZ,
  \alpha\in\pi\}\\
  X_\pi^+&:=&\{\la\in\mh^*\mid \langle \la,\check{\alpha}\rangle\in\mZ_{\geq
  0}, \alpha\in\pi\}.
\end{eqnarray*}
Recall that there is a natural bijection
\begin{eqnarray}
\label{simples}
  \{\text{finite dimensional irreducible $\cU(\p)$-modules}\}&\oneone&
  X_\pi^+,
\end{eqnarray}
by mapping a module to its highest weight. Let $\tilde{E}(\la)$ be
the irreducible module corresponding to $\la\in X_\pi^+$. For any
$\pi$-deformation algebra $T'$ define the {\it $T'$-deformed
  (generalised) Verma module} with highest weight $\la\in X_\pi^+$ as
\begin{eqnarray*}
  M_{T'}^\p(\la):=\cU(\mg)\otimes_{\cU(\p)}(\tilde{E}(\la)\otimes T').
\end{eqnarray*}
This is a $\cU(\mg)\otimes T'$-module, where $T'$ is just acting on $T'$ by
multiplication. Given a $\cU(\mg)\otimes T'$-module $M$ and $\la\in\mh^*$, we
denote by $$M_{T'}^\la=\{m\in M\mid h.m=\varphi(\la(h)+h)m,
\forall\;h\in\mh\}$$ the {\it
  $\la$-weight space} of $M$. (Here $\la(h)+h$ is considered as an element of
  $T$ via the map $i'\circ i\circ p$ and $\varphi(\la(h)+h)$ is an element of
  $T'$.)

\subsection{The deformed (parabolic) category $\cO$}
Let $\cO_{T}^\p$ be the full subcategory of the category of
$\cU(\mg)\otimes T$-modules defined by the set of objects $M$
satisfying
\begin{itemize}
\item $M$ is finitely generated,
\item $M=\bigoplus_{\la\in\mh^*} M_{T}^\la$ as $T$-module, and
\item $(\cU(\p)\otimes T)m$ is a finitely generated $T$-module for all $m\in
  M$.
\end{itemize}
In particular, $M_{T}^\p(\la)\in\cO_{T}^\p$ for any $\la\in
X_\pi^+$. Note that the third condition is equivalent to saying
that $M$ is locally $\cU(\mg_\pi\oplus\mn^\pi)\otimes T$-finite.

If we replace $T$ by $\mC$ in all the definitions then $\cO_{T}^\p$ is the ordinary parabolic category $\cO^\p$ as defined in \cite[Section 3]{RC}. In particular, $\cO^\mb_\mC$ is the ordinary BGG-category $\cO^\mb$ from \cite{BGG}. The generalised Verma modules $M^\p_\mC(\la)$ are abbreviated as $M^\p(\la)$.

\subsection{Weight and root decompositions of $\cO_T^\p$}

If $\la\in\mh^*$ we denote by $\overline{\la}$ its class in $\mh^*/ X$ and by
$\tilde{\la}$ its class in $\mh^*/\mZ R$. We have the following {\it weight
decomposition}: $  \cO_{T}^\p=\bigoplus_{\Lambda\in\mh^* /X
}\cO^\p_{T,\Lambda},$ where $M\in\cO_{T,\Lambda}^\p$ if and only if
$M_{T}^\la\not=\{0\}\Rightarrow \overline \la=\Lambda,$ and the finer {\it root
decomposition}
%\begin{eqnarray*}
$  \cO_{T}^\p=\bigoplus_{\Lambda\in\mh^*/\mZ R}\cO^\p_{T,\Lambda},$
%\end{eqnarray*}
where $M\in\cO_{T,\Lambda}^\p$ if and only if $M_{T}^\la\not=\{0\}\Rightarrow
\tilde\la=\Lambda.$ Both decompositions are far away from giving block
decompositions, hence we have to refine them once more.

\subsection{Central character and block decomposition}
Let $\xi^\sharp:\mathfrak{Z}\rightarrow S$ be the Harish-Chandra homomorphism
normalised by $\xi^\sharp(z)-z\in\cU(\mg)\mn$, where $\mn=\mn^{\emptyset}$. Hence
$\xi^\sharp: \mathfrak{Z}\rightarrow S^{(W, \cdot)}$, the $W$-invariants under the
dot-action. This is the comorphism for $\xi:\mh^*\rightarrow
\operatorname{Max}(\cZ)$ which induces a bijection between $(W,\cdot)$-orbits
of $\mh^*$ and maximal ideals $\operatorname{Max}(\cZ)$ of $\mathfrak{Z}$. For
any weight $\la$, we have $+\la:\mh^*\rightarrow\mh^*$, $\mu\mapsto \la+\mu$.
Let $\la^\sharp:S\rightarrow S$ be the corresponding comorphism.

Since deformed Verma modules are generated by their highest weight space we
 have a canonical isomorphism of rings $\END_{\mg\otimes
 T}\big(M_{T}^\p(\la)\big)=T$. Let $\chi_\la:\mathfrak{Z}\otimes T\rightarrow T$ be such that
 $z.m=\chi_\la(z)m$ for any $z\in\mathfrak{Z}\otimes T$, $m\in
 M_{T}^\p(\la)$.
Explicitly, the morphism $\chi_\la$ is given by
 \begin{eqnarray}
\label{chila}
   \mathfrak{Z}\otimes T\stackrel{(\la^\sharp\circ\xi^\sharp)\otimes id}\longrightarrow S\otimes
   T\stackrel{i^\sharp\otimes id}\longrightarrow S(\mh^\pi)\otimes
   T\stackrel{j\otimes\test{id}}\longrightarrow T\otimes T\stackrel{m}\longrightarrow T,
 \end{eqnarray}
where $i:(\mh^\pi)^*\rightarrow \mh^*$ is the canonical
embedding, $j$ is the canonical embedding into its localisation,
and $m$ is the multiplication map (for details see \cite[Section 2]{Sperv}).

Let $\mu\in \mh^*$ and consider the support $\op{supp}
M_{T}^\p(\la)\subseteq\op{Spec}(\mathfrak{Z}\otimes T)$ of $M_{T}^\p(\la)$
viewed as a $\mathfrak{Z}\otimes T$-module. By definition, $\op{supp}
M_{T}^\p(\la)$ is the closed subset of the spectrum of $\mathfrak{Z}\otimes T$
given by all prime ideals containing $\op{Ann}_{\mathfrak{Z}\otimes T}
(M_{T}^\p(\la))=\op{Ann}_{\mathfrak{Z}\otimes T}
(M_{T}^\p(\la)^\la)=\KER\chi_\la$. Unlike in the non-deformed situation,
$\KER\chi_\la$ does not need to be a maximal ideal. However, the homomorphism
theorem implies that there is a homeomorphism between the spectrum of $T$ and
$\op{supp} M_{T}^\p(\la)$. Hence, $\op{supp} M_{T}^\p(\la)$ contains exactly
one closed point, namely $\xi(\la)\otimes T+\mathfrak{Z}\otimes \mm$, since $T$
is local with maximal ideal $\mm$. If $\la$, $\mu\in X_\pi^+$ then
\begin{eqnarray*}
 \op{supp} M_{T}^\p(\la)\cap \op{supp}
 M_{T}^\p(\mu)\not=\emptyset\Leftrightarrow \xi(\la)=\xi(\mu)\Leftrightarrow
 \la\in W\cdot\mu.
\end{eqnarray*}

For $\chi$, a maximal ideal of $\mathfrak{Z}$, let $\cO^\p_{T,\chi}$ be the
full subcategory of $\cO_{T}^\p$ given by all objects having support contained
in $\cap_{\xi(\mu)=\chi}\op{supp} M_{T}^\p(\mu)$. For $\la\in\mh^*_{dom}$ let
$\cO^\p_{T,\la}=\cO^\p_{T,\xi(\la)}\cap\cO^\p_{T, \Lambda}$, such that
$\Lambda=\tilde{\la}$. We have the following {\it `block' decomposition}:
\begin{eqnarray*}
  \cO_{T}^\p&=&\bigoplus_{\la\in\mh^*_{dom}}\cO^\p_{T,\la}
=\bigoplus_{\la\in\mh^*_{dom}\cap X_\pi^+}\cO^\p_{T,\la}.
\end{eqnarray*}

Strictly speaking, this is not a block decomposition, since the summands might
decompose further. This is, however, not the case if $\p=\mb$ or $\la=0$, where
$\cO^\p_{T,\la}$ is in fact a block. Since we are mainly interested in this
case we call it `block' decomposition.

\subsection{The ordinary parabolic category $\cO_0^\p$}
\label{ordinary}
Let us stop for a moment and recall the structure of the principal
block $\cO_0^\p$ of $\cO_{\mC}^\p=\cO^\p$ from \cite{RC}. The generalised Verma
modules in $\cO_0^\p$ are exactly the $M^\p_\mC(\la)$, where $\la\in
X_\pi^+\cap W\cdot0$, or in other words $\la$ is of the form
$\la=w\cdot0$, where $w\in W^\p$. The simple objects in $\cO_0^\p$
are exactly the simple quotients $L(w\cdot0)$, $w\in W^\p$ of
these generalised Verma modules. (There is only one finite
dimensional simple module, namely the trivial module.) The
category $\cO_0^\p$ has enough projectives; for $w\in W^\p$ let
$P^\p(w\cdot0)$ be the projective cover of $L(w\cdot0)$ in
$\cO^\p_0$.

\subsection{Specialisations:  $\mC$ and $\cQ=\op{Quot}{T}$}
\label{specialisation}

For any morphism $f: T\rightarrow T'$ of deformation algebras let $\cO^\p_{T',\la}$ denote the image of $\cO^\p_{T,\la}$ under the specialisation functor
$_-\otimes_TT'$.

From \cite[Proposition 2.6 and Section 2.4]{CT} it follows that if $f:
T\rightarrow T/\mm=\mC$ is the canonical projection then the image of
$\cO_T^\p$ under the $\mC$-specialisation is the ordinary parabolic category
$\cO^\p_\mC$ with the usual decomposition
$\cO_\mC^\p=\bigoplus_{\la\in\mh^*_{dom}\cap X_\pi^+}\cO^\p_{\mC,\la}$ (see
e.g. \cite[4.4]{Ja2}). It follows directly from the definitions that
$M_T^\p(\la)\otimes_T\mC\cong M^\p_\mC(\la)$ as $\mg$-modules.

On the other hand we could consider the specialisation functor
$_-\otimes_T\cQ$, where $\cQ=\op{Quot}{T}$ is the quotient field of $T$. We
identify
$(\mh\otimes_\mC\cQ)^*:=\HOM_\cQ(\mh\otimes\cQ,\cQ)=\mh^*\otimes_\mC\cQ$. Let
$\tau\in(\mh\otimes\cQ)^*$ be the tautological weight, i.e. restricted to $\mh$
it is just the projection onto $\mh^\pi\subseteq S(\mh^\pi)\subset
T\subset\cQ$. From the definitions we have $M^\p_T(\la)\otimes_T\cQ\cong
M_\cQ^\p(\la+\tau)$ as $\mg\otimes\cQ$-modules. If $\la\in X_\pi^+$ for all
$\alpha\in\pi$ then
$\langle\la+\tau,\check\alpha\rangle=\langle\la,\check{\alpha}\rangle\in\mZ$ by
definition of $\mh^\pi$ and $\tau$, hence $\la+\tau$ is an admissible weight
for the Lie algebra $\mg\otimes\cQ$. If $\alpha\in\triangle-\pi$ then
$\tau(\alpha)\not=0$ (since the elements from $\triangle$ are linearly
independent) and therefore $\langle \la+\tau,\check\alpha\rangle\notin\mZ$ for
any $\alpha\in\triangle-\pi$. It follows that  $M_\cQ^\p(\la+\tau)$ is simple
(\cite[(1.17)]{Ja2}). In particular, the image of $\cO_T^\p$ under the
$\cQ$-specialisation functor is semisimple with simple objects
$M_\cQ^\p(\la+\tau)$, $\la\in X_\pi^+$.

\subsection{Translation functors}
Let $\la$, $\mu\in\mh^*_{dom}$ such that $\mu-\la\in X(R)$. For any deformation algebra $T'$ let
\begin{eqnarray*}
  \theta_{\la,T'}^\mu:\cO^\p_{T',\la}\rightarrow \cO^\p_{T',\mu}
\end{eqnarray*}
be the translation functor defined as $M\mapsto\op{pr}_\mu(M\otimes E)$, where
$E$ is the finite dimensional $\mg$-module with extremal weight $\mu-\la$ and
$\op{pr}_\mu$ is the projection to the summand $\cO^\p_{T',\mu}$. Since
$\cO_{T'}^\p$, the direct sum of all blocks, is closed under tensoring with
finite-dimensional $\mg$-modules, the definition makes sense. Obviously,
$\theta_{\la,T'}^\mu$ commutes with base change, i.e, there is a natural
isomorphism
\begin{eqnarray*}
  \theta_{\la,T'}^\mu(M\otimes_T T')&\cong&(\theta_{\la,T}^\mu M)\otimes_T T'.
\end{eqnarray*}

\subsection{Deformed projectives}
For the reader's convenience we recall some fundamental properties of the deformed parabolic categories, but omit the proofs. The arguments can be found in \cite{Sperv}, \cite{SHC}, \cite{CT}.

\begin{prop}
Let $T'$ be any $\pi$-deformation algebra. \label{general}
  \begin{enumerate}
  \item \label{1} The category $\cO_{T'}^\p$ has enough projectives.
  \item \label{2} The category $\cO_{T'}^\p$ is closed under taking direct
  summands and finite direct sums.
   \item \label{4} If $\la\in X_\pi^+\cap\mh_{dom}^*$ then $M_{T'}^\p(\la)$ is
  projective in $\cO_{T',\la}^\p$ and in  $\cO_{T'}^\p$.
   \item \label{5} Any projective module in  $\cO_{T'}^\p$ is obtained by
  applying translation functors to some $M_{T'}^\p(\la)$, $\la\in
  X_\pi^+\cap\mh^*_{dom}$, taking finite direct sums and direct summands.
   \item \label{9} Any projective object in $\cO_{T'}^\p$ has a {\it
  Verma-flag}, i.e. a filtration with subquotients isomorphic to various
  deformed generalised Verma modules.
  \item \label{3} The weight spaces of projective objects in $\cO_{T'}^\p$ are
  free $T'$-modules of finite rank.
  \item \label{8} The specialisation functor ${}_-\otimes_T\mC$ defines a
  bijection between the (indecomposable) projective objects in $\cO_{T}^\p$
  and the (indecomposable) projective objects in $\cO_{\mC}^\p$.
  \item \label{6}
If $M$, $N\in\cO_{T'}^\p$ are projective, then $\HOM_{\mg\otimes T'}(M,N)$ is a free $T'$-module of finite rank.
  \item \label{7} If $M$ $N\in\cO_T^\p$ are projective then the canonical map
    \begin{eqnarray*}
      \Psi:\HOM_{\mg\otimes T}(M,N)\otimes_T T'&\cong&\HOM_{\mg\otimes
      T'}(M\otimes_T T',N\otimes_TT')\\
      f\otimes t&\mapsto&f\otimes t\cdot\test{id}.
    \end{eqnarray*}
  is an isomorphism.
  \end{enumerate}
\end{prop}

If $P\in\cO_{T'}^\p$ has a Verma flag then we denote by
$[P:M_T^\p(\la)]$ the mul\-tipli\-ci\-ty of a $M_T^\p(\la)$ as a subquotient of an (arbitrarily chosen) Verma flag of $P$.
Note that this number is stable under changes of the deformation
ring (by Proposition~\ref{general}~\eqref{3}). For simplicity we
will restrict our attention to the principal blocks
$\cO_{T',0}^\p$ in the following.

\subsection{Commutativity of the endomorphism rings}
In this section we will use the deformation theory to obtain (as
an easy application) the commutativity of certain endomorphism
rings:

\begin{prop}
\label{comm}
  Let $P\in \cO_{0}^\p$ be an indecomposable projective module.
  \begin{enumerate}
  \item \label{first} If $[P :M^\p(\la)]\leq 1$ for any $\la\in\mh^*$ then $\END_\mg(P)$ is commutative.
  \item \label{second} If  $\p=\mb$ then the following are equivalent:
  \begin{enumerate}[(i)]
  \item  $[P :M^\mb(\la)]\leq 1$ for any $\la$,
  \item  $\END_\mg(P)$ is commutative,
  \item $\mathfrak{Z}$ surjects onto $\END_\mg(P)$ canonically.
  \end{enumerate}
  \end{enumerate}
\end{prop}

\begin{proof}
  Let $P_T$ be the $T$-deformation of $P$ given by
  Proposition~\ref{general}~\eqref{8}. From Proposition~\ref{general}~\eqref{7}
  we get an isomorphism of rings $$\END_{\mg\otimes T}(P_T)\otimes_T T'\cong\END_{\mg\otimes
  T}(P_T\otimes_T T')$$ for any $T$-algebra $T'$. If we choose $T'=\mC$ then
  the commutativity of  $\END_{\mg\otimes T}(P_T)$ implies the commutativity of $\END_\mg(P)$. On the other hand, we could choose $T'=\cQ$, the ring
of fractions of $T$. The category $\cO_{0,\cQ}$ is semisimple, with simple
objects being the $\cQ$-specialised deformed Verma modules
(Section~\ref{specialisation}). They all have commutative endomorphism rings
isomorphic to $\cQ$. Set $J=\{\la\in\mh^*\mid [P_T:M_T^\p(\la)]\not=0\}$. By
our assumption on the multiplicities we get $P_T\otimes_T \cQ\cong
\oplus_{\la\in J}M_T^\p(\la)\otimes_T \cQ\cong
  \oplus_{\la\in J} M_{\cQ}^\p(\la+\tau)$, and
  $\END_{\mg\otimes\cQ}(P\otimes_T{\cQ})\cong
\oplus_{\la\in J}\cQ$ is commutative. Proposition~\ref{general}~\eqref{6}, \eqref{7} provide an
inclusion
  \begin{eqnarray*}
    \END(P_T)&\longrightarrow&\END(P_T)\otimes_T\cQ\cong\END(P\otimes_T\cQ)\\
f&\longmapsto& f\otimes 1.
  \end{eqnarray*}
  So, $\END_{\mg\otimes T}(P_T)$ is a subring of a
  commutative ring, hence itself commutative. The first part of the proposition follows. The second is \cite[Theorem 7.1]{Stquiv}.
\end{proof}

We don't know if Proposition~\ref{comm} \eqref{second} is true for general
$\p$. A famous example for a module $P$ satisfying the conditions of the
proposition in case $\p=\mb$ is the `antidominant projective module'
$P(w_0\cdot0)\in\cO_0^\mb$ (see Section~\ref{centreacts} below). Here are
further examples (see also Proposition~\ref{Irving}):

 \begin{prop}
\label{maxparab}
   Let $\mg=\mathfrak{gl}_n$ and $\p_\pi$ a maximal parabolic,
   i.e. $\pi=\triangle-\{\alpha_s\}$ for some simple reflection $s$. Then
   $\END_\mg(P)$ is commutative for any indecomposable projective object
   $P\in\cO_0^\p$.
 \end{prop}

 \begin{proof}
   With this choice of a parabolic subalgebra, the assumptions of
   Proposition~\ref{comm}~\eqref{first} are satisfied (\cite[Theorem 5.1]{Brenti}).
 \end{proof}

In Section~\ref{TL} we will see that the endomorphism rings appearing in Proposition~\ref{maxparab} are of the form $(\mC[X]/(X^2))^{\otimes k}$ for some $k\in\mZ_{\geq 0}$.

\section{The centre as a module for the Weyl group}
\label{centreacts}

Recall, that the centre $\cZ(\cA)$ of an abelian category
$\cA$ is the endomorphism ring of the identity functor
$id=id_\cA$. If, for instance, $\cA\cong\MOD-A$, the
category of finitely generated right modules over some
unitary $\mC$-algebra $A$ then the centre of $\cA$ is
naturally isomorphic to the (ordinary) centre $\cZ(A)$ of
the algebra $A$. The isomorphism associates to a natural
transformation $f$ the value $f_A(1)$.

\subsection{The cohomology ring of the flag variety}

Let us consider for a moment the case $\p=\mb$. It is well-known that the
centre of $\cO^\mb_{\mC,0}$ is naturally isomorphic to
$\mathfrak{C}=S/(S^W_+)$, the ring of coinvariants (\cite[Endomorphismensatz
and Struktursatz]{Sperv} together with \cite[Theorem 5.2 (2)]{MSSerre}). The
natural action of the Weyl group $W$ on $S=S(\mh)$ gives rise to an action of
$W$ on $\mathfrak{C}$. In the deformed situation, the picture is similar: the
centre of the deformed category $\cO^\mb_{0,T}$ is naturally isomorphic to
$S\otimes_{S^W} T$ (\cite[Theorem 9, Corollary 1]{SHC}), and hence carries
obviously the structure of a $W$-module. To get an explicit description of the
isomorphism we first note that each element of the centre of $\cO^\mb_{0,T}$
defines an element of the endomorphism ring $E$ of the `antidominant
projective' in $\cO^\mb_{0,T}$ by restriction, defining an isomorphism between
the centre of the category and $E$ (see e.g. \cite[Theorem 1.8]{StTQFT}). On
the other hand, Soergel showed in \cite[Theorem 9]{SHC} that $E$ is canonically
isomorphic to $T\otimes_{T^W} T=S\otimes_{S^W} T$. Moreover, since the
specialisation functor ${}_-\otimes_T\mC$ maps $E$ surjectively onto the
endomorphism ring of the antidominant projective module in $\cO^\mb_{\mC,0}$,
the principal block of the ordinary category $\cO$, the centre of
$\cO^\mb_{T,0}$ maps surjectively onto the centre of $\cO^\mb_{\mC,0}$.

Let $\cB=G/B$ be the flag variety corresponding to $\mg$ (i.e. the variety of
Borel subalgebras in $\mg$) and $H^*(\cB)$ its cohomology algebra with complex
coefficients. The Weyl group acts on $H^*(\cB)$. Note that $\mathfrak{C}$ has
an even $\mZ$-grading coming from the grading on $S$, where $\mh$ is
concentrated in degree two. We recall the following well-known fact (see e.g.
\cite[Section 4.1]{CoPr} or \cite[Proposition 7.2]{Springertrig}):

\begin{prop}
\label{centreSoerg}
  There is a $W$-equivariant isomorphism of graded algebras $\psi:\mathfrak{C}\cong H^*(G/B)$, and $\mathfrak{C}\cong \mC[W]$ as $W$-modules.
\end{prop}

Via the natural isomorphism $\mathfrak{C}=\cZ(\cO_0^\mb)$, the centre of
$\cO_0^\mb$ inherits an action of $W$ giving rise to the regular
representation. In the following we will explain how this $W$-action on the
centre can be obtained via braid group actions on the bounded derived category
$\cD^b(\cO_0^\mb)$
 - inducing a $W$-action on the centre of $\cO^\mb_{\mC,0}$ and
 then finally also on the centre of $\cO^\p_{\mC,0}$.

\subsection{Braid group actions on the centre of a category}

Before we pass to derived categories, we want to give the main idea behind this
braid group action on the centre of $\cO_0^\p$ by first assuming a simplified
situation. Let $\cC$ be an abelian $\mC$-linear category. Let $F:\cC\rightarrow
\cC$ be a functor. Assume $F$ is invertible. Then the centre of $\cC$ is
isomorphic to $\END(F)$ in two ways: First by mapping an element $c$ in the
centre to $F(c)$ and secondly by mapping $c$ (naively) to the endomorphism
given by multiplication with $c$. Now given $c$ in the centre of $\cC$ there is
a unique $c'$ in the centre of $\cC$ such that $F(c)$ is given by
multiplication with $c'$. In particular, there is an automorphism
$\Psi_F:\cZ(\cC)\rightarrow \cZ(\cC)$ which maps $c$ to $c'$. Hence $c$ and $c'$ are related by the formula $c'_{F(M)}=F(c_M)$ for any object
$M$.

Assume now $G$ is a group acting on $\cC$, i.e. for any $g\in G$, there is an
(invertible) functor $F_g:\cC\rightarrow\cC$ such that $F_e\cong\test{id}$,
$F_{gh}\cong F_g\circ F_h$. We get the corresponding automorphisms $\Psi_{F_g}$
which give rise to an action of $G$ on the centre of $\cC$. For more detail we
refer to \cite{KhoSpringer}.

\subsection{The Irving shuffling functors}
\label{shuffle} Let $s\in W$ be a simple reflection and choose $\la\in\mh^*$ an
integral weight with stabiliser $\{e,s\}$. Let
$\theta_s=\theta_\la^0\theta_0^\la:\cO^\mb_0\rightarrow\cO^\mb_0$ be the
translation functor through the $s$-wall. Let $a_s:\ID\rightarrow\theta_s$ be
the adjunction morphism. Consider the functor $\mathrm{C}_s=\COKER(a_s)$. This
is a right exact functor such that its left derived functor $\cL \mathrm{C}_s$
induces an equivalence on the bounded derived category $\cD^b(\cO_0)$
(\cite[Theorem 5.7]{MS}). It is quite easy to see that they satisfy braid
relations in the weak sense (\cite[Theorem 2]{KM} and \cite[Section 6.5]{MOS}),
which means if we have a braid relation $st\ldots =ts\ldots$ then there is an
isomorphism of functors $\cL \mathrm{C}_s\cL \mathrm{C}_t\ldots\cong\cL
\mathrm{C}_t\cL \mathrm{C}_s\ldots$. (Although this weak version is enough for
our purposes we want to point out that Rouquier showed that the isomorphisms of
functors can be chosen in a compatible way, \cite{Rouquier}.)

Since the translation functors preserve the parabolic categories, the functors
$\mathrm{C}_s$ induce functors $\mathrm{C}_s:\cO^\p_0\rightarrow \cO^\p_0$ and
the left derived functors $\cL \mathrm{C}_s$ are auto-equivalences of
$\cD^b(\cO^\p_0)$ (cf. \cite[Section 4]{MS}).

Each of the categories $\cO^\p_0$ is equivalent to $\MOD-A^\p$ for some finite
dimensional algebra $A^\p$ (see e.g. \cite[Section 2.1]{Stquiv}). Under this
equivalence, the functors $C_s$ become so-called tilting functors, given by
tensoring with some tilting complex (see \cite{Rickard}, \cite[Section 5]{MS}).
Hence we have a braid group action via tilting auto-equivalences on
$\cD^b(\cO^\p_0)$. Since these equivalences are given by tilting complexes we
obtain an induced braid group action on the centre of the underlying abelian
category $\cO^\p_0$ (\cite[Theorem 9.2]{Rickardcentre}).

\subsection{The action of the Weyl group on $\cZ(\cO^\mb_0)$}
Let us now construct this action explicitly. We first consider the case $\p=\mb$ and recall some results from \cite{Sperv}:
Let $P(w_0\cdot0)\in\cO^\mb_0$ be the projective cover of the simple
Verma module $M(w_0\cdot0)=L(w_0\cdot0)$. Consider Soergel's
{\it Strukturfunktor}
$$\mV=\HOM_\mg(P(w_0\cdot0),_-):\cO^\mb_0\longrightarrow
\MOD-\END_\mg(P(w_0\cdot0)).$$

By Soergel's Endomorphismensatz we
have $\END_\mg(P(w_0\cdot0))\cong \mathfrak{C}=\mathfrak{C}^{op}$
canonically and under this identification we get a functor $\mV:\cO^\mb_0\rightarrow\mathfrak{C}-\MOD$. There is an isomorphism
$\mV\theta_s\cong\Theta_s\mV$, where $\Theta_s:\mathfrak{C}-\MOD\rightarrow\mathfrak{C}-\MOD$, $M\mapsto
\mathfrak{C}\otimes_{\mathfrak{C}^s}M$  and $\mathfrak{C}^s$ denotes the
$s$-invariants of $\mathfrak{C}$. Under this isomorphism, the  adjunction morphism $a_s:\test{id}\rightarrow \theta_s$ corresponds to the morphism $a_s$ given by
\begin{eqnarray}
\label{comult}
a_s(M):\quad\mathfrak{C}\otimes_\mathfrak{C}M&\rightarrow&
\mathfrak{C}\otimes_{\mathfrak{C}^s}M,\nonumber\\
1\otimes m&\mapsto&X\otimes
m+1\otimes Xm,\quad m\in M,\quad X=\check{\alpha}_s
\end{eqnarray}
(\cite[Lemma 8.2]{StDuke}).
Let $\COKER_s$ be the functor of taking the
cokernel of $a_s:\op{id}\rightarrow\Theta_s$.
By construction $\mV\mathrm{C}_s\cong \COKER_s\mV$ when restricted to the additive subcategory of $\cO^\mb_0$ generated by $P(w_0\cdot0)$.

\begin{lemma}
\label{action} Let $c\in \mathfrak{C}$ with the corresponding element $m_c$ of
the centre of $\mathfrak{C}-\MOD$ given by multiplication with $c$. Then
$\COKER_s(m_c) = m_{s(c)}$.
\end{lemma}

\begin{proof}
It is enough to check this on the regular module $\mathfrak{C}$. Let
$\mathfrak{C}\stackrel{a_s}\rightarrow \mathfrak{C}\otimes_{\mathfrak{C}^s}
\mathfrak{C}\stackrel{p}\rightarrow \COKER_s(\mathfrak{C})$ be the defining
sequence of $\COKER_s(\mathfrak{C})$. Note that $\COKER_s(\mathfrak{C})\cong
\mathfrak{C}$ as $\mathfrak{C}$-module, and is generated by the image of
$1\otimes 1$ under $p$. Let $c\in\mathfrak{C}$, homogeneous of degree one. Then
we have
\begin{eqnarray}
\label{ccc} \Theta_s(m_c)(1\otimes 1)=1\otimes m_c(1)=1\otimes c.
\end{eqnarray}
On the other hand $s(c)-c=r X$ for some $r\in\mC$ and $X=\check{\alpha}_s$.
Since $\mathfrak{C}$ is a free $\mathfrak{C}^s$ module on basis $1$, $X$, we
therefore get
\begin{eqnarray*}
1\otimes c=\frac{1}{2}\big(1\otimes(s(c)+c)-1\otimes
(s(c)-c)\big)=\frac{1}{2}\big((s(c)+c)\otimes 1-1\otimes r X\big),
\end{eqnarray*}
because $c+s(c)$ is $s$-invariant. However,
\begin{eqnarray*}
\frac{1}{2}\big((s(c)+c)\otimes 1-1\otimes r X\big)
&\equiv&\frac{1}{2}\big((s(c)+c)\otimes 1+r X\otimes 1\big)\\
&=&\frac{1}{2}\big((s(c)+c+s(c)-c)\otimes 1\big)\\
&=&s(c)\otimes 1,
\end{eqnarray*}
where $\equiv$ means equality modulo the image of $a_s$. The lemma follows.
\end{proof}

Let us summarise: On $\cZ(\cO_0^\p)$, the centre of $\cO_0^\p$, there is an
action of the braid group $B_W$ which underlies $W$. This action is induced
from the braid group action of the left derived functors of Irving's shuffling
functors (Section~\ref{shuffle}).

\begin{theorem}
\label{centreb} Let $\p\subseteq\mg$ be any parabolic subalgebra containing
$\mb$.
  \begin{enumerate}
    \item The action of the braid group on $\cZ(\cO_0^\mb)$ factors through $W$.
    \item The canonical isomorphism $\cZ(\cO_0^\mb)=\mathfrak{C}$ is $W$-equivariant.     \item The canonical restriction morphism $\cZ(\cO_0^\mb)\rightarrow\cZ(\cO_0^\p)$ is $B_W$-equivariant. In particular, the image becomes a $W$-module.
  \end{enumerate}
\end{theorem}

\begin{proof}
  The first two statements hold because of Lemma~\ref{action} and the natural
  identification of the centre with the endomorphism ring of the antidominant
  projective module by restriction (\cite[Theorem 5.2 (2)]{MSSerre}).
  The last statement follows directly from the definition of the braid group
actions.
\end{proof}

\begin{remark}
\label{semisimple}
  {\rm Theorem~\ref{centreb} (1) and (2) hold analogously for the deformed categories $\cO_{T,0}^{\p}$. If we consider the corresponding semisimple category $\cO_{\cQ,0}^\p$ then we have isomorphisms
      \begin{eqnarray*}
        \cZ(\cO_{\cQ,0}^\p)&\cong&\bigoplus_{x\in W^\p}\END_{\mg\otimes\cQ}M_\cQ^\p(x\cdot0+\tau)=\bigoplus_{w\in W^\p}\cQ\\
        z&\mapsto& (z_x)_{x\in W^\p},
      \end{eqnarray*}
      where $z_x\in \cQ$ is the image of the natural transformation $z$ applied to $M_\cQ(x\cdot0+\tau)$ evaluated at $1\otimes 1\otimes 1\otimes 1\in M_\cQ(x\cdot0+\tau)$. The $\cQ$-version of $\cL\mathrm{C_s}$ maps $M_\cQ(x\cdot0)$ to $M_\cQ(xs\cdot0)$ if $x$, $xs\in W^\p$ and to $M_\cQ(x\cdot0)[1]$ otherwise. Hence $s(z_x)=z_{xs}$ if $xs\in W^\p$ and $s(z_x)=z_x$ otherwise. Therefore, $\cZ(\cO_{\cQ,0}^\p)\cong \mC[W]\otimes_{\mC[W_\p]}\mC_{triv}$ as $W$-module.
}
\end{remark}

\section{Type $A$: the centres and the Springer fibres}
\label{typeA}

From now on we will stick to the special case where
$\mg=\mathfrak{gl}_n$ with standard Borel subalgebra $\mb$ given
by the upper triangular matrices. We would like to generalise Proposition~\ref{centreSoerg} to the parabolic case.

Let $G=GL(n,\mC)$ with Borel subgroup $B$ given by the upper triangular
matrices and Lie algebra $\mb$. Then $\cB=G/B$ is the variety of complete flags
in $\mC^n$. Let $\p\subseteq\mg$ be a parabolic subalgebra containing $\mb$,
hence $\p$ is given by upper diagonal $(\mu_1, \mu_2,\dots, \mu_r)$-block
matrices, where $\mu=(\mu_1,\mu_2,\dots,\mu_r)$ is a composition of $n$. For
example, $\p=\mb$ corresponds to the composition $(1^n)$ of $n$. Any maximal
parabolic subalgebra corresponds to a two part composition. Let
$\la(\mu)=\la(\p)$ be the partition obtained from $\mu$ by reordering the
parts. Let $x=x_\p\in G$ be the nilpotent element of $G$ in Jordan Normal Form
such that the Jordan blocks are of size $\mu_1$, $\mu_2$, etc.  Let $\cB_\p$ be
the Springer fibre corresponding to $x$, that means the subvariety of $\cB$ of
all flags fixed by the unipotent element $u=\operatorname{Id}+x$. Note that
$\cB_\mb=\cB$.

\subsection{The main result: centres via cohomology}
Let $H^*(\cB_\p)$ denote the cohomology of $\cB_\p$. Springer defined  an $S_n$-action on $H^*(\cB_\p)$ and proved that the top part $H^*(\cB_\p)_{top}$ is the irreducible representation of $S_n$ corresponding to the partition $\la(\p)$ (\cite{Springertrig} or \cite{HottaSpringer}). The embedding of $\cB_\p$ into
$\cB$ induces a morphism $h:H^*(\cB_\mb)\rightarrow H^*(\cB_\p)$,
which is surjective (\cite{HottaSpringer}) and $W$-equivariant (\cite{CoPr}, \cite{Tanisaki}).

\begin{theorem}
\label{main} Let $\mg=\mathfrak{gl}_n$ with parabolic subalgebra
$\p$ containing $\mb$ and Weyl group $W=S_n$.

\begin{enumerate}
\item \label{mainone}
The canonical map $H^*(\cB)=\mathfrak{C}=\cZ(\cO_0)\rightarrow
  \cZ(\cO^\p_0)$ factors through $H^*(\cB_\p)$ and induces an isomorphism of rings
  \begin{eqnarray*}
    \Phi_\p: H^*(\cB_\p)\iso  \cZ(\cO^\p_0).
  \end{eqnarray*}
\item\label{maintwo}
$\Phi_\p$ is $W$-equivariant and $ \cZ(\cO^\p_0)\cong\mC[W]\otimes_{\mC[W_\p]}\mC_{triv}$ as $W$-module.
\item\label{mainfour}  Up to isomorphism, $H^*(\cB_\p)$ and $\cZ(\cO^\p_0)$ only depend on $\la(\p)$.

\end{enumerate}
\end{theorem}

\begin{proof}
Let $R$ denote the regular functions on $\mh^*\oplus U$, where
  $U:=\{\la\in\mh^*\mid \la(\mh_\pi)=0\}=(\mh^\pi)^*$. Hence $R=S\otimes
  S(\mh^\pi).$ We fix the standard basis $\epsilon_i$, $1\leq i\leq n$ for $\mh^*$ with its set of fundamental weights $\omega_i=\sum_{k=1}^i\epsilon_i$. The $\omega_{i_1}, \ldots, \omega_{i_r}$ contained in $U$ form a basis of $U$. Let
  $$I=\{f\in R\mid f\big(w(\la),\la\big)=0, \text{ for all $\la\in U, w\in
  W$}\}\subseteq R\subseteq S\otimes T$$ and put $K=\{f(_-,0)\mid f\in I\}\subset S=S(\mh)$.\\

{\it Claim 1:}  $K=\KER\big(S\rightarrow \mathfrak{C}\rightarrow
H^*(\cB_{\p})\big)=:\KER$.\\

We start by showing that $K$ contains $\KER$. Thanks to \cite{Tanisaki} we have
an explicit set of generators for $\operatorname{ker}$: Let
$\mu'=(\mu_1',\mu_2',\ldots, \mu_{r'}')$ be the dual partition of $\la(\mu)$.
If we identify $S$ with $\mC[x_1,\ldots x_n]$ in the usual way by taking the
dual standard basis vectors $\epsilon_i^*$ as generators, then $\operatorname{ker}$ is generated by
all $l$-th elementary symmetric functions $e_l(\mathcal{X})$, $k>0$, $l>0$
where
$$\mathcal{X}\subseteq\{x_1,\ldots, x_n\},\quad |\mathcal{X}|=k,\quad k\geq
l>k-(\mu_{n-k+1}'+\mu_{n-k+2}'+\ldots+\mu_{r'}').$$ Therefore, it is enough to
show that these $e_l(\mathcal{X})$ are contained in $K$. Thanks to the
$W$-invariance, we only have to consider the cases where $\mathcal{X}$ consists
of the first $k$ variables $x_1, x_2,\ldots, x_k$. Taking the dual basis of the $\omega_{i_j}\in U$, $1\leq j\leq r$, we will identify
$R=\mC[x_1,x_2,\ldots, x_n]\otimes\mC[y_1,y_2,\ldots, y_r]$, and for any choice
of $k$ and $l$ from the allowed range construct a polynomial $f=f_{k,l}\in R$
with the following properties
\begin{enumerate}
\item $f(x_1,x_2,\ldots, x_n,0,0,\ldots, 0)=e_l(x_1,x_2,\ldots, x_k),$
\item $f(b_1,b_2,\cdots, b_n,a_1,a_2,\ldots, a_r)=0$ for any point $(b_1, b_2,
\ldots, b_n)$ where $\mu_i$ of the coordinates are equal to $a_i$ for $1\leq
i\leq r$.
\end{enumerate}
Then $f\in I$ and so $e_l(x_1,x_2,\ldots, x_k)=f(x_1,x_2,\ldots,
x_n,0,0,\ldots, 0)\in K$.\\

Let now $k$ and $l$ be fixed from the allowed range. The construction of this
polynomial goes along the lines of \cite{GP}. For $1\leq i\leq r$ let $m_i$ be
the maximum of $\mu_i-n+k$ and zero, and define $d:=\sum_{i=1}^r m_i$. Note
that $d=(\mu_{n-k+1}'+\mu_{n-k+2}'+\ldots+\mu_{r'}')$ as above. To construct
the polynomial $f$ we first define $P(t)=\prod_{i=1}^{k}(t+x_i),$ a polynomial
in $t$ with coefficients from $\mC[x_1,x_2,\ldots x_k]$, and consider the
polynomial $Q(t)=\prod_{i=1}^{r}(t+y_i)^{m_i}$ of degree $d$ in $t$ with
coefficients from $\mC[y_1,y_2,\ldots, y_r]$. We do the long division of $P(t)$
by $Q(t)$ and get $P(t)=q(t)Q(t)+r(t),$ where $r(t)=\sum_{s=0}^{d-1}r_s\;t^s$
is a polynomial in $t$ with coefficients being homogeneous polynomials $r_s$ in
the $x$'s and $y$'s. We claim that
$f=r_{k-l}$ does the job.

If we set all $y_i=0$, then $Q(t)=t^d$, and $r_s=e_{k-s}(x_1, \ldots, x_k)$,
hence $r_{k-l}=e_{l}(x_1, \ldots x_k)$ and $(1)$ follows. Now let $(a_1,
\ldots, a_r)\in \mC^r=U$ and ${\bf b}=(b_1, b_2, \ldots b_n)\in \mC^n=\mh^*$
such that $\mu_i$ of the coordinates are equal to $a_i$ for $1\leq i\leq r$,
hence at least $m_i=\mu_i-(n-k)$ of the first $k$ coordinates are equal to
$a_i$. In particular, $P(t)$ evaluated at the first $k$ coordinates of ${\bf
b}$, is divisible by $Q(t)$ if $y_i=a_i$, i.e. if $Q(t)$ is evaluated at ${\bf
b}$. Hence, $r(t)$ becomes zero when evaluated at $x_i=b_i$, $y_i=a_i$, and so
$(2)$ holds. This implies $\operatorname{ker}\subseteq K$.

It is left to show that the inclusion is in fact an equality. Let $u\in U$ be a
generic point and define $K_u=\{f(_-,u)=0\mid f \in I\}\subseteq S$. Since $u$
is generic, $S/K_u$ is the coordinate algebra of $|W/W_\mathfrak{p}|$ distinct
points in $\mh^*$, and hence $\operatorname{dim} (S/{K_u})=|W/W_\mathfrak{p}|$,
and also $\operatorname{dim} (S/gr K_u)=|W/W_p|$, where the associated graded
$gr$ is taken with respect to the canonical grading of $S$. Moreover, one can
easily see that $K\subseteq gr(K_u)$ by mapping $f(_-,0)$ to $f(_-,u)$.
Altogether, the natural surjection from $S$ to $S/K$ factors through
$S/\operatorname{ker}$ and the natural surjection from $S$ to $S/gr K_u$
factors through $S/K$. Therefore, $\operatorname{dim} (S/\KER)\geq
\operatorname{dim} (S/K)\geq \operatorname{dim} (S/gr K_u)=|W/W_\mathfrak{p}|$.
By the main result of \cite{Tanisaki}, we have $\operatorname{dim}
(S/\KER)=|W/W_\mathfrak{p}|$,
hence all the dimensions agree and Claim 1 follows.\\

{\it Claim 2:} $K\subseteq\KER\big(S\rightarrow\cZ(\cO_0)\rightarrow
  \cZ(\cO^\p_0)\big)$.\\

Let $P_T$ be a projective generator of $\cO_{T,0}^\p$. We consider the
following commutative diagram:
\begin{eqnarray*}
\xymatrix { \mathfrak{Z}\otimes T \ar[r]^{\xi^\sharp} \ar[d]^{\alpha} &S^W\otimes T
\ar[rr]^{h} \ar[d]^{\beta} &&\END_{\mg\otimes T}(P_T)
\ar[d]^\gamma\\
\mathfrak{Z}\otimes \mathcal{Q} \ar[r]^{\xi^\sharp} &S^W\otimes \mathcal{Q}
\ar[rr]^{h} &&\END_{\mg\otimes\mathcal{Q}}(P_{\mathcal{Q}})}
\end{eqnarray*}
The two maps denoted $h$ are given by applying the product of the two factors to the
module (so that $h\circ\xi^\sharp$ is the canonical map), whereas the vertical maps are the obvious inclusions (using
Proposition~\ref{general}~\eqref{7}). By Section~\ref{specialisation} we have
an isomorphism
$$\epsilon:\quad\END_{\mg\otimes\mathcal{Q}}(P_T\otimes_T{\mathcal{Q}})\cong\END_{\mg\otimes{\mathcal{Q}}}\left(\bigoplus
a_x M^\p_\cQ(x\cdot0+\tau)\right),$$ where $a_x$ is the multiplicity of
$M^\p_T(x\cdot0)$ as a subquotient of a Verma flag of $P_T$ (see Section~\ref{specialisation}). The map
$\gamma$ is injective (by Proposition~\ref{general}~\eqref{6}). Hence the
kernel of the upper {\it can} is the kernel of $can\circ\beta$.

An element $z\otimes t\in\mathfrak{Z}\otimes T$ acts on
$M^\p_\cQ(x^{-1}\cdot0+\tau)$ by multiplication with
$(x^{-1}\cdot0+\tau)(\xi^\sharp(z))\tau(t)$. On the other hand
$(x^{-1}\cdot0+\tau)\circ\xi^\sharp(z) =\xi^\sharp(z)(x^{-1}\cdot0+\tau)= (x\cdot
x^{-1}\cdot0+x(\tau))\circ \xi^\sharp(z) =x(\tau)\circ\xi^\sharp(z)$, so $z\otimes t$ acts by multiplication with $x(\tau)\circ\xi^\sharp(z)$. On the other hand
\begin{equation}
\label{action2}
j_x:\quad S\otimes T\stackrel{x^\sharp\otimes id}\longrightarrow S\otimes T\stackrel{p_\pi\otimes id}\longrightarrow S(\mh^\pi)\otimes T\stackrel{mult}\longrightarrow T\longrightarrow\cQ,
\end{equation}
where $x^\sharp:S\rightarrow S$ is the comorphism given by the action of $x$.
Hence $z\otimes t$ acts by multiplication with $\tau\circ j_x(\xi^\sharp(z)\otimes t)$.

To pass from $\mathfrak{Z}\otimes T$ to $S\otimes T$ note that the map $\xi^\sharp$ induces an isomorphism $\psi:\mathfrak{Z}\otimes T/{{\xi^\sharp}^{-1}(J)} \cong S\otimes T/{J}$, where $J=\cap_{x\in W}\KER j_x\subseteq S\otimes T$, \cite[p.429]{Sperv}. (Note that we use here the assumption that our block $\cO_0^\p$ is regular, and so $\xi$ is etale at $0$.)

Let now $f\in K$ and find $\tilde{f}\in I$ such that $f=\tilde{f}(_-,0)$. Using the explicit formula before \eqref{action2} and the map $\psi$ it follows that $\tilde{f}$ acts on $M_\cQ(x^{-1}\cdot0+\tau)$ by multiplication with the
function $\la\mapsto\tilde{f}(x(\la),\la)$, hence by zero. So, $\tilde{f}$ is
in the annihilator of $P_T$. By definition of $K$, we can write
$\tilde{f}=f\otimes 1+g\in S\otimes T$, where $g\in S\otimes\mm$. So, $g$ induces an
endomorphism on $P_T$ which specialises to the zero endomorphism of $P_\mC$.
Claim 2 follows.\\

Therefore, the canonical map $H^*(\cB)=\mathfrak{C}=\cZ(\cO_0)\rightarrow
  \cZ(\cO^\p_0)$ factors through $H^*(\cB_\p)$ and induces the
  map $\Phi_\p$. The latter is surjective by \cite[Theorem 5.11]{Brundan}.
  The injectivity will be proved at the end of the section. Let us assume for the moment we proved this already (so \eqref{mainone}
  holds).
  Thanks to Theorem~\ref{centreb}, $\Phi_\p$ is $W$-equivariant, and  Statement~\eqref{maintwo} is true if $\p=\mb$ by Proposition~\ref{centreSoerg}. Hence the image of $\Phi_\p$, i.e. $\cZ(\cO^\p_0)$,  is isomorphic to $\mC[W]\otimes_{\mC[W_\p]}\mC_{triv}$ as $W$-module by
  \cite[Corollary 8.5]{HottaShim} or \cite[Theorem 1]{Tanisaki}.
  Statement~\eqref{maintwo} of the theorem follows.
  Part~\eqref{mainfour} is clear from \cite{Tanisaki}, but also has a direct
  proof from the categorical side by \cite[Theorem 5.4 and Theorem 5.2]{MSSerre}.
\end{proof}

In the following two subsections we will prove the injectivity of $\Phi_\p$.
This result follows in fact directly by dimension arguments, since Brundan
showed that the dimension of the centre of $\cO^\p_0$ equals the dimension of
$H^*(\cB_\p)$ (\cite{Brundan}). However, our approach gives a distinguished
basis of the top degree part of $H^*(\cB_\p)$ and shows in a nice way how the
categorification of the symmetric group action comes into play.

\subsection{A generalised antidominant projective module}

The missing part in the proof of Theorem~\ref{main} will be deduced from several non-trivial results which we recall first. We start with the following fact:

\begin{prop}
\label{Irving}
  Let $\mg=\mathfrak{
gl}_n$ and $\p$ be any parabolic subalgebra. Then there is
  always an indecomposable projective module $P\in\cO_0^\p$ such that
  \begin{itemize}
  \item $P$ is injective, and
  \item the natural action of $\mathfrak{C}$ defines a surjection $\mathfrak{C}\surj\END_\mg(P)$, in particular $\END_\mg(P)$ is commutative.
  \end{itemize}
 \end{prop}

 \begin{proof}
   By \cite[Corollary p.327]{IS} there is an integral dominant weight $\la$ such that
   $\cO_\la^\p$ contains a simple projective module $N=L(\mu)$,
   $\mu=x\cdot\la$ for some $x\in W$. Hence $N$ is
   also injective and a Verma module. Then $\theta_\la^0(N)$ is projective and
   injective, and has a (generalised) Verma flag satisfying the assumptions of
   Proposition~\ref{comm} (\cite[4.13 (1)]{Ja2}). Hence there is some
   projective and injective module $P=\theta_\la^0N$ with commutative endomorphism
   ring. If
   $\Gamma:\cO_\la\rightarrow \cO_\la^\p$ is the functor which picks out the
   largest quotient contained in $\cO_\la^\p$ then $\Gamma M(\mu)=N$ . Since
   $\Gamma$ commutes with translation functors, we have $\Gamma\theta_\la^0M(\mu)\cong
   \theta_\la^0\Gamma M(\mu)\cong\theta_\la^0 N=P$. On the other hand, $\Gamma$
   commutes (by definition) with the action of the centre and thanks to the
   existence of the canonical
   projection $\END_\mg(\theta_\la^0M(\mu))\rightarrow\END(\theta_\la^0N)$, it
   is enough to show that $\mathfrak{C}$ surjects onto $\END_\mg(\theta_\la^0M(\mu))$
   naturally. However, $M(\mu)\cong T_xM(\la)$, where $T_x$ is the twisting
   functor as studied in \cite{AL}, \cite{AS}. Now, $T_x$ commutes with the
   action of the centre (see the definition of the functors in \cite{AS}), it is therefore enough to show that $\mathfrak{C}$
   surjects onto $\END_\mg(\theta_\la^0M(\la))$ naturally. However,
   $\theta_\la^0M(\la)$ satisfies the assumptions of Proposition~\ref{comm}
   for $\p=\mb$. Then the statement follows from Proposition~\ref{comm}~\eqref{second}.
 \end{proof}

 \begin{definition}
An object $M\in\cO_0^\p$ is called {\it projective-injective} if it is both, projective and injective in $\cO_0^\p$.
\end{definition}
Proposition~\ref{Irving} ensures the existence of projective-injective modules.

\subsection{Grading and Loewy length}
\label{gradings}

Let $P\in\cO_0^\p$ be projective, then $\END_\mg(P)$ has a natural non-negative $\mZ$-grading induced from the Koszul grading (\cite{BGS}) of $A^\p=\END_\mg(P_{gen})$, where $P_{gen}\in\cO_0^\p$ is a minimal projective generator. This Koszul grading induces a $\mZ$-grading on $\cZ(\cO_0^\p)$, the centre of $\cO_0^\p$.

\begin{prop}
\label{ll}
    Assume $\mg$ is any reductive complex Lie algebra and $\p$ some
    parabolic subalgebra. Let $P_i$, $i\in I$,  be a complete system of representatives for the
    isomorphism classes of indecomposable projective-injectives in $\cO_0^\p$.
    \begin{enumerate}
    \item \label{ll1} The centre of $\cO_0^\p$ is the
    centre of $\END_\mg(\oplus_{i\in I} P_i)$.
    \item The Loewy lengths of all $P_i$ agree and equal the maximal possible
    Loewy length $\op{ll}$. The maximal degree of $\END_\mg(P_i)$, considered as a graded ring, is equal to $\op{ll}-1$.
    \item Consider $\cZ^\p:=\cZ(\cO_0^\p)$ as a graded ring with its top
    degree part $\cZ^\p_{top}$. Then $\DIM_\mC\cZ^\p_{top}=|I|$ and $top=\op{ll}-1$.
    \end{enumerate}
 \end{prop}

 \begin{proof}
The first two statements of the Proposition are proved in \cite[Theorem 5.2]{MSSerre}. Now consider $D:=\END_\mg(\oplus_{i\in I} P_i)$ as a graded ring.
Let $\tilde{P}_i$ be a graded lift of $P_i$ (in the sense of \cite{BGS} or \cite[Section 3.1]{StGrad}).
Then the radical, socle and the grading filtrations of $\tilde{P}_i$ agree up to a shift of the grading (\cite[Proposition 2.4.1]{BGS}), since $P_i$ has simple head and simple socle (the latter by \cite[Appendix]{Irself}). If now $f\in D$ is of maximal degree then $f$ is contained in the span of the maps $g_i$ where $g_i$ maps the head of $P_i$ to its socle and is zero on all other summands by the second statement. On the other hand, the $g_i$ are all contained in the centre of $B$. Therefore we have an isomorphism of vector spaces $D_{top}\cong\cZ^\p_{top}$ and the proposition follows.
\end{proof}

\begin{lemma}
\label{top}
Let $\mg=\mathfrak{gl}_n$ and $\p$ be any parabolic subalgebra. Then the top degree of $\cZ(\cO_0^\p)$ agrees with the top degree of $\cH^*(\cB_\p)$.
\end{lemma}

\begin{proof}
From \cite[Proposition and Corollary 3.1]{IS} (see also Remark~\ref{thanks}) we have an explicit formula for the Loewy length of a projective-injective module in $\cO_0^\p$, hence for the top degree of $\cZ(\cO_0^\p)$ by Proposition~\ref{ll}. The formula agrees with \cite[Lemma 1.3]{HottaShim} and implies the assertion.
\end{proof}

\begin{lemma}
\label{numerology}
  Let $\mg=\mathfrak{gl}_n$ and $\p$ be any parabolic subalgebra. Then the
  following numbers coincide
  \begin{itemize}
  \item the number of isomorphism classes of indecomposable
  pro\-jec\-tive-in\-jec\-tive modules in $\cO_0^\p$,
\item the dimension of $\cZ^\p_{top}$, where $\cZ^\p$ is the centre of  $\cO_0^\p$,
  considered as a graded algebra,
 \item the number of irreducible components of $\cB_\p$, and hence the dimension of $H^*(\cB_\p)_{top}$.
\item the dimension of the irreducible representation $S^{\la(\p)}$ of the symmetric group $S_n$.
  \end{itemize}
\end{lemma}

\begin{proof}
  The first two agree by Proposition~\ref{ll}, the last two by Springer's construction of irreducible $S_n$-modules (\cite[Proposition 7.1]{HottaShim}). Thanks to the main result of \cite{Irself}, the indecomposable projective-injective modules in $\cO^\p_0$ are indexed by elements of some right cell of $S_n$. The statement follows, since the cell modules are exactly the irreducible $S_n$-modules (\cite{Na}).
\end{proof}

The importance of the bijections in the previous lemma becomes apparent in the fact that the category of indecomposable projective-injective modules in  $\cO_0^\p$ together with the translation functors $\theta_s$, $s\in W$ a simple
  reflection, categorifies the irreducible representation $S^{\la(\p)'}$ of $S_n$, where $\la(\p)$ denotes the dual partition of $\la(\p)$.

 More precisely: Consider the additive category $\cC^\p$ of projective-injective modules in $\cO_0^\p$.  By Lemma~\ref{numerology}, the complexified (split) Grothendieck group $K_0(\cC^\p)$ (that is the complexification of the free abelian group generated by the isomorphism classes $[M]$ of objects $M$ in $\cC^\p$ modulo the relation $[M]+[N]=[M\oplus N]$) is isomorphic to the corresponding Specht module as a complex vector space. Moreover, $\cC^\p$ is stable under translations $\theta_s$ through walls. Since the functors $\theta_s$ are exact, they induce endomorphisms $[\theta_s]$ of $K_0(\cC^\p)$. Let $T_s=[\theta_s]-id:K_0(\cC^\p)\rightarrow K_0(\cC^\p)$ then the statement is the following:

\begin{prop}[\cite{KMS}]
\label{categorif}
Let $\mg=\mathfrak{gl}_n$ with $\mb$ the standard Borel and $\p$ a parabolic subalgebra. Let $\la(\p)'$ be the dual partition of $\la(\p)$ and $S^{\la(\p)'}$ the corresponding Specht module. Then there is an isomorphism of right $S_n$-modules
\begin{eqnarray*}
  \epsilon:\quad S^{\la(\p)'}\cong K_0(\cC^\p),
\end{eqnarray*}
The $S_n$-module structure on the right hand side is induced by the $T_s$'s.
\end{prop}

Let us finish the proof of Theorem~\ref{main} as follows:

\begin{theorem}
\label{injectivity}
  Let $\mg=\mathfrak{gl}_n$ and let $\p$ be any parabolic subalgebra.
  Then the map $\Phi_\p$ from Theorem~\ref{main} is an inclusion which induces an isomorphism of $S_n$-modules $H^*(\cB_\p)_{top}\cong\cZ^\p_{top}$.
\end{theorem}

\begin{proof}
  We know that $\Phi_\p$ is a homomorphism of $\mathfrak{C}$-modules. Hence, it is enough to show that $\Phi_\p$ is injective when restricted to the socle of $H:=H^*(\cB_\p)$. By \cite[Theorem 6.6 (vi)]{IainBaby} the socle of $H$ agrees with the part of highest degree $H_{top}$. We have $\Phi_\p(H_{top})\subseteq\cZ^\p_{top}$ by Lemma~\ref{top} and get $\Phi_\p(H_{top})\not={0}$ by Proposition~\ref{Irving} and Proposition~\ref{ll}. On the other hand $\Phi_\p$ is $B_W$-equivariant and even $W$-equivariant onto its image (Theorem~\ref{centreb}). Since $H_{top}$ is an irreducible $W$-module, $\Phi_\p$ defines an inclusion  $H^*(\cB_\p)_{top}\inj\cZ^\p_{top}$, and $\Phi_\p$ is injective. Moreover, $H^*(\cB_\p)_{top}\inj\cZ^\p_{top}$ must be an isomorphism by Lemma~\ref{numerology}.
\end{proof}

As a consequence (independent of \cite[Theorem 5.11]{Brundan})) we
get the following categorical construction of the Springer representations
\begin{cor}
\label{SpringerIrr}
There is an isomorphism $\cZ^\p_{top}\cong S^{\la(\p)}$ of $S_n$-modules.
\end{cor}

\begin{proof}
  This follows directly from Theorem~\ref{injectivity} and Springer's construction of the irreducible $S_n$-modules (\cite{Springer}), since $\Phi_\p$ is $W$-equivariant.
\end{proof}

\subsection{A few remarks on the singular case}

Let still $\mg=\mathfrak{gl}_n$.
Let $\nu\in\mh^*$ be an integral dominant weight. Let $W_\nu=\{w\in W\mid w\cdot\nu=\nu\}$.
Then $\cZ(\cO_\nu^\mb)=\mathfrak{C}^{W_\nu}$, the $W_\nu$ invariants of $\mathfrak{C}$ (\cite{Sperv}).
Let $\p^\nu$ be the parabolic subalgebra of $\mg$ such that $W_{\p^\nu}=W_\nu$ and denote by
$P^\nu$ the subgroup of $GL(n,\mC)$ with Lie algebra $\p^\nu$. Let $\cP^\nu=G/P^\nu$ be the corresponding partial flag variety. The canonical projection $\cB\surj \cP^\nu$ gives rise to an inclusion $i^\nu_\mb: H^*(\cP^\nu)\hookrightarrow H^*(\cB)$ with image $\mathfrak{C}^{W_\nu}$ (\cite[Lemma 8.1]{HottaShim}). For any parabolic $\p$ of $\mg$ let $\cP^\nu_\p$ be the associated fixed point variety. The image of the corresponding inclusion $$i^\nu_\p: H^*(\cP^\nu_\p)\hookrightarrow H^*(\cB_\p)$$ are the $W_{\nu}$-invariants $H^*(\cB_\p)^{W_{\nu}}$ of $H^*(\cB_\p)$ (\cite[Lemma 8.1]{HottaShim}).\\

Theorem~\ref{main} and Lemma~\ref{numerology} generalise as follows:
\begin{theorem}
\label{theosingular}
Let $\mg=\mathfrak{gl}_n$ with the standard Borel $\mb$ contained in a fixed parabolic $\p=\p_{\pi}$. Let $\nu\in\mh^*$ be a dominant integral weight.
\begin{enumerate}
\item \label{theosingularone}
The canonical map
$H^*(\cP^\nu)=\mathfrak{C}^{W_{\nu}}=\cZ(\cO^\mb_\nu)\rightarrow
  \cZ(\cO^\p_\nu)$ factors through $H^*(\cP^\nu_\p)$ and induces a ring homomorphism
  \begin{eqnarray*}
    \Phi_\p^\nu: H^*(\cP^\nu_\p)\longrightarrow\cZ(\cO^\p_\nu)
  \end{eqnarray*}
\item\label{theosingularthree}  $\Phi_\p^\nu$ is a surjection and an isomorphism if $\p=\mb$.
\end{enumerate}
\end{theorem}

\begin{proof}[Proof of Theorem~\ref{theosingular}]
The canonical projection of $\cB_\p$ onto $\cP^\nu_\p$ induces an inclusion
$H^*(\cP^\nu_\p)\rightarrow H^*(\cB_\p)$. By \cite[Lemma 8.1]{HottaShim} we
know that the image are exactly the $W_{\nu}$-invariants of $H^*(\cB_\p)$.
Assume $z$ is in the kernel of the canonical map $H^*(\cP^\nu)   \rightarrow
H^*(\cP^\nu_\p)$. Assume $z$ acts non-trivially on $\cO_\nu^\p$. Then there is
some module $M\in\cO_\nu^\p$, such that $z$ acts non-trivially on $M$. Hence
$\theta_0^\nu(z)\in\END_\mg(\theta_0^\nu M)$ is non-trivial, since
$\theta_0^\nu$ is exact and does not annihilate modules. By \cite[Lemma
8]{Sperv} $\theta_0^\nu(z)$ is just given by multiplication with
$z\in\mathfrak{C}$. Hence $z$ acts non-trivially on $\cO_0^\p$. This is a
contradiction to Theorem~\ref{main} and the first statement of the theorem
follows.

The map $\Phi_\mb^\nu$ is an isomorphism by \cite[Endomorphismensatz]{Sperv}.
On the other hand $\Phi_\p^\nu$ is surjective by \cite[Theorem 5.11]{Brundan}.
\end{proof}

\begin{lemma}
\label{projinj} Let $\mg$ be any reductive complex Lie algebra with Borel $\mb$
and some parabolic subalgebra $\p\supset\mb$. Let $\nu\in\mh^*$ be a dominant
integral weight. Then $P\in\cO^\p_\nu$ is indecomposable projective-injective
if and only if so is $\theta_\nu^0P\in\cO^\p_0$.
\end{lemma}

\begin{proof}
Let $L\in\cO_0^\p$ be a simple module. Then $\theta_0^\nu L=L'$ is simple or
zero (\cite[4.12 (3)]{Ja2}), and $\HOM_\mg(\theta_\nu^0
P,L)=\HOM_\mg(P,\theta^\nu_0L)=\HOM_\mg(P,L')\not=0$ only if $P$ is the
projective cover of $L'$. Using again (\cite[4.12 (3)]{Ja2}) we deduce that
$\theta_\nu^0 P$ has simple top, and is therefore indecomposable. Since
$\theta_\nu^0$ does not annihilate any module, $P$ is indecomposable if and
only if $\theta_\nu^0P$ is indecomposable. Since $\theta_0^\nu\theta_\nu^0$ is
isomorphic to a direct sum of copies of the identity functor
(\cite[4.13(2)]{Ja2} and \cite[Theorems 3.3 and 3.5]{BG}) and translation
functors preserve projectivity and injectivity,  $P$ is projective-injective if
and only if so is $\theta_\nu^0P$.
\end{proof}

\begin{remark}{\rm
\label{thanks} Proposition~\ref{categorif} together with
Proposition~\ref{projinj} give the dimension of the top degree part of
$\cZ(\cO_\nu^\p)$, namely the number of standard $\nu$-tableaux of shape
$\la(\p)'$, the dual partition of $\la(p)$. Using a graded version of
$\theta_\nu^0$ (in the sense of \cite[Definition 3.3]{StGrad}) one can deduce
from  Proposition~\ref{projinj} a formula for the top degree $top$ of
$\cZ(\cO_\nu^\p)$: Consider $\la(\p)'$ and $w_0^{\p'}$ the longest element in
the corresponding symmetric group $S_{\la(\p)'_1}\times \ldots \times
S_{\la(\p)'_r}$. If $\nu$ is regular then $top=2(l(w_0^{\p'}))$ by
\cite[Proposition and Corollary 3.1]{IS}. If $\nu$ is not necessarily regular,
let $w_0^\nu$ be the longest element in $W_\nu$. Then
$top=2(l(w_0^{\p'})-l(w_0^\nu))$, because a graded version of $\theta_\nu^0$
adds $2l(w_0^\nu)$ degrees if we apply it to a simple module (see \cite[Theorem
8.2 (4)]{StGrad} for a special case). If either $\p\not=\mb$ or $\nu$ not
regular then the number of simple objects in $\cO_\nu^\p$ is in fact (see
\cite[Theorem 2]{Brundan}) strictly smaller than the dimension of
$H^*(\cP^\nu_\p)$ (Section~\ref{ordinary} and \cite[Remark 8.6]{HottaShim}). In
particular, the map $\Phi_\p^\nu$ is not injective in these cases
(\cite[Theorem 2]{Brundan}). A detailed description of these resulting proper
quotients of $H^*(\cP^\nu_\p)$ can be found in \cite{Brundan2}. }
\end{remark}

\subsection{The maximal parabolic case}
\label{max}

The proof of the injectivity of $\Phi_\p$ in Theorem~\ref{main} was quite involved. In case $\mg=\mathfrak{gl}_n$ and $\p\not=\mg$ is a maximal parabolic subalgebra, there is an alternative proof which we want to present now.

\begin{lemma}
\label{symmetric}
  Let $B$ be a finite dimensional complex algebra. Assume $B$ is symmetric, i.e. there is a non-degenerate associative symmetric $\mC$-bilinear form $\mathbf{b}:B\times B\rightarrow \mC$. Then there is an isomorphism of vector spaces
  \begin{eqnarray*}
    \cZ(B)&\cong&(B/[B,B])^*\\
    z&\mapsto&\mathbf{b}(z,_-).
  \end{eqnarray*}
\end{lemma}

\begin{proof}
  Let $z\in \cZ(B)$, the centre of $B$ and $a, b\in B$. Then $\mathbf{b}(z,ab-ba)=\mathbf{b}(z,ab)-\mathbf{b}(z,ba)=\mathbf{b}(z,ab)-\mathbf{b}(zb,a)=\mathbf{b}(z,ab)-\mathbf{b}(a,zb)=\mathbf{b}(z,ab)-\mathbf{b}(a,bz)=\mathbf{b}(z,ab)-\mathbf{b}(ab,z)=0$. Hence $\mathbf{b}(z,_-)\in (B/[B,B])^*$. On the other hand if $\mathbf{b}(z,_-)\in (B/[B,B])^*$ then  $\mathbf{b}(z,ab-ba)=0$ for any $a,b\in B$. Hence $\mathbf{b}(za,b)=\mathbf{b}(z,ab)=\mathbf{b}(z, ba)=\mathbf{b}(ba,z)=\mathbf{b}(b,az)=\mathbf{b}(az,b)$ for all $b\in B$ and so $az=za$, since $\mathbf{b}$ is non-degenerate, and therefore $z\in\cZ(B)$. The claim of the lemma follows.
\end{proof}

\begin{theorem}
\label{conjA}
   Let $\mg=\mathfrak{gl}_n$ and $\p_\pi$ a maximal parabolic,
   i.e. $\pi=\triangle-\{\alpha_s\}$. Then $\Phi_\p$ is an isomorphism.
\end{theorem}

\begin{proof}
The surjectivity is given by \cite{Brundan}. Let $P_i$, $i\in I$,  be a complete system of representatives for the
    isomorphism classes of indecomposable projective-injectives in $\cO_0^\p$. Set $P=\oplus_{i\in I}P_i$ and let $P_T$ be the $T$-deformation of $P$ given by Proposition~\ref{general}~\eqref{8}. Then $A_T=\END_{\mg\otimes T}(P_T)$ is a free $T$-module of finite rank (Proposition~\ref{general}~\eqref{6}). Let $D_T:=[A_T, A_T]$. Since $\p=\p_{\pi}$ is assumed to be maximal parabolic, $\mh^\pi$ is one-dimensional and $S(\mh^\pi)$ is a principal ideal domain. Therefore, $D_T$ is a free $T$-module as well.
We have canonical isomorphisms $A_T\otimes_TT'\cong \END_{\mg\otimes T}(P_T)\otimes_TT'\cong \END_{\mg\otimes T'}(P_T\otimes_T T')$ for $T'=\mC$ or $T'=\cQ$ (Proposition~\ref{general}~\eqref{7}). Set $A_T\otimes_TT'=A_{T'}$ for $T'=\mC$ or $T'=\cQ$ and note that $D_{T'}:=D_T\otimes_T {T'}=[A_T, A_T]\otimes_T {T'}$ surjects onto $[A_{T'}, A_{T'}]$ canonically. For $T=\cQ$ we even have an isomorphism, since $D_T$ is free as $T$-module. We deduce that
    \begin{eqnarray}
\label{dims}
      &\DIM_\mC\END_\mg(P)=\DIM_\cQ\END(P_T\otimes_T\cQ)=\op{rank}_T\END(P_T)=\op{rank}_T A_T&\nonumber\\
 &\DIM_\mC D_\mC\leq\DIM_\cQ D_\cQ=\op{rank}_T D_T.&
    \end{eqnarray}
Since $\END_{\mg\otimes\cQ}(P_\cQ)$ is a product of $|W^\p|$ matrix rings (see Section~\ref{specialisation}), we have $\DIM_\cQ(A_\cQ)-\DIM_\cQ D_\cQ=\DIM_\cQ(A_\cQ/D_\cQ)=|W^\p|$. Hence the equations~\eqref{dims} imply $|W^\p|\leq\DIM_\mC(A_\mC)-\DIM_\mC D_\mC=\DIM_\mC(A_\mC/D_\mC)$. Since the algebra $A_\mC$ is symmetric (\cite[Theorem 5.2]{MSSerre}), we get $$\DIM_\mC\cZ(A_\mC)\geq |W^\p|.$$ Since the map $\Phi_\p$ is surjective, $\cZ(A_\mC)$ is a quotient of $ \mC[W]\otimes_{\mC[W_\p]}\mC_{triv}$ as $W$-module by Theorem~\ref{centreb} and \cite[Corollary 8.5]{HottaShim}. Comparing the dimensions we are done.
\end{proof}

\section{Diagrammatic approach and Khovanov homology}
\label{TL}

In this section we prove Theorem~\ref{conjTQFT} and Theorem~\ref{Theothree} from the introduction. We proceed as follows: First we recall the definition of the algebras $\cH^n$. Then we state some combinatorial results which will be used to give a purely diagrammatic description of $\cO^\p(\mathfrak{gl}_n)$ for maximal parabolic subalgebras $\p\not=\mg$. We explicitly describe how Braden's presentations can be transformed into this diagrammatic framework, where almost everything is computable. This will finally improve the presentation of \cite{Braden} drastically in several ways: we are able to see the Koszul grading and obtain a usual Ext-quiver with {\it homogeneous} relations, we give a very easy recipe to compute dimensions of homomorphism spaces between projective modules, and deduce that the endomorphism rings of projective-injective modules are all isomorphic.

For the rest of the paper we fix $n\in\mZ_{>0}$ and $\mg=\mathfrak{gl}_{2n}$
with standard Borel $\mb$ and $\p=\p_n$ the parabolic subalgebra where $W_\p=S_n\times S_n$ and denote $\cO_0^{\p}(\mathfrak{gl}_{2n})$ by
$\cO_0^{n,n}$. This is the category which plays an important role.

\subsection{Khovanov's algebras $\cH^n$}

We recall the basic definitions from \cite{Khovanov}, but refer to that paper
for details. From now on let $R:=\mC[X]/(X^2)$ be the ring of dual numbers. This is a commutative Frobenius algebra, hence defines a 2-dimensional topological quantum field theory $\cF$. In other words, $\cF$ is a monoidal functor from the category of oriented cobordisms between $1$-manifolds to the category of finite dimensional complex vector spaces. The Frobenius algebra structure of $R$ is given by
\begin{itemize}
\item the associative multiplication ${\bf m}:R\otimes R\rightarrow R$, $r\otimes s\mapsto rs$,
\item the comultiplication map $\Delta: R\rightarrow R\otimes R$, $1\mapsto X\otimes 1+1\otimes X$, $X\mapsto X\otimes X$ (note that this is just a special case of \eqref{comult} for $\mg=\mathfrak{gl}_2$),
\item the unit map $\epsilon: \mC\rightarrow R$, $1\rightarrow 1$,
\item the counit or trace map $\delta: R\rightarrow\mC$, $1\mapsto 0$, $X\mapsto 1$.
\end{itemize}

The functor $\cF$ associates to $k$ disjoint circles the vector space
$R^{\otimes k}$, to the cobordisms of `pair of pants shape' the multiplication
map ${\bf m}$ and the comultiplication $\Delta$ respectively. To the cobordisms
connecting one circle with the empty manifold, $\cF$ associates the trace map
$\delta$ and the unit map $\epsilon$.
\begin{figure}[htb]
  \centering
  \includegraphics{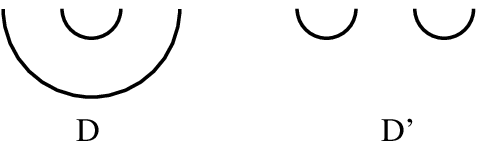}
  \caption{The elements of $\op{Cup}(2)$}
  \label{fig:cups}
\end{figure}

\begin{figure}[htb]
  \centering
  \includegraphics{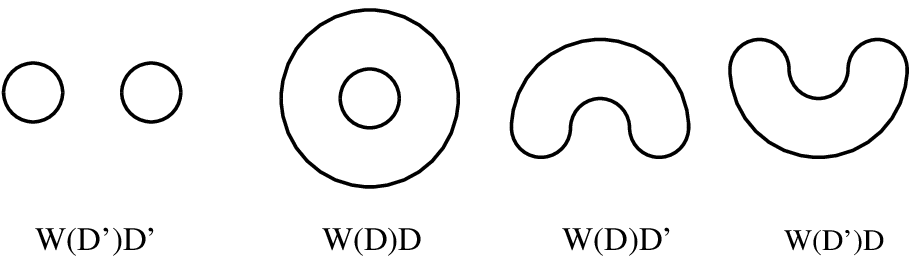}
  \caption{Gluing two cup diagrams gives a collection of circles}
  \label{fig:gluing}
\end{figure}

\begin{figure}[htb]
  \centering
  \includegraphics[scale=0.9]{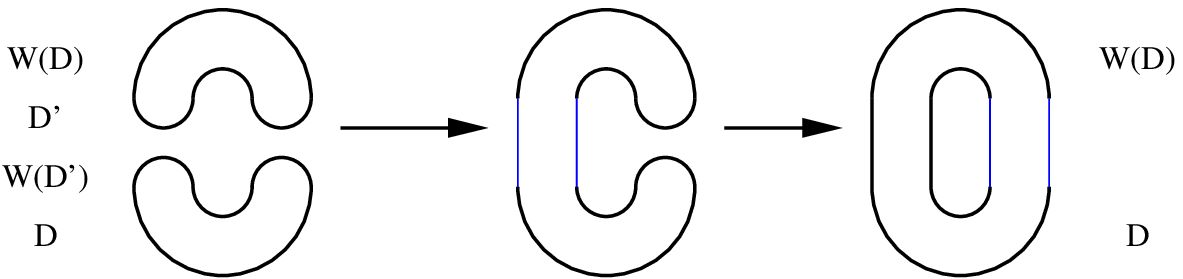}
  \caption{Composition of maps: $R\otimes R\stackrel{\op{m}}\longrightarrow R\stackrel{\Delta}{\longrightarrow} R\otimes R$}
  \label{fig:contract}
\end{figure}

\begin{figure}[htb]
  \centering
  \includegraphics[scale=0.9]{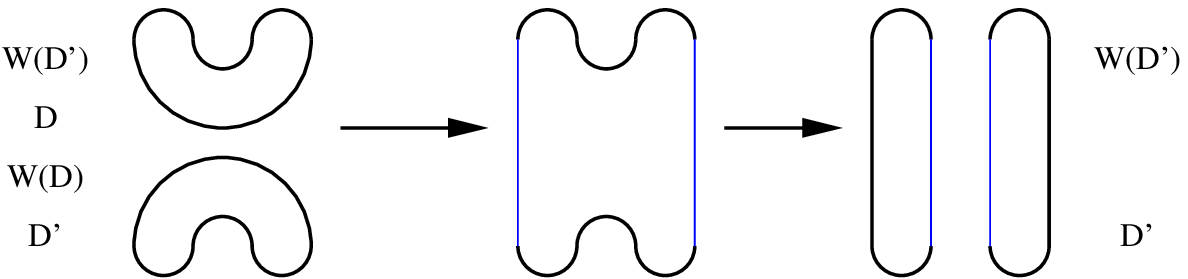}
  \caption{Composition of maps: $R\otimes R\stackrel{\op{m}}\longrightarrow R\stackrel{\Delta}{\longrightarrow} R\otimes R$}
  \label{fig:contractb}
\end{figure}

Let $\op{Cup}(n)$ be the set of crossingless matchings of $2n$ points (see Figure~\ref{fig:cups}). For $a$, $b\in\op{Cup}(n)$ let $W(b)$ be the reflection of $b$ in the horizontal axis and $W(b)a$ the closed 1-manifold obtained by gluing $W(b)$ and $a$ along their boundaries (see Figure~\ref{fig:gluing}).
Given $a,b,c\in\op{Cup}(n)$, there is the cobordism from $W(c)bW(b)a$ to $W(c)a$ which contracts $bW(b)$ (see Figures~\ref{fig:contract} and \ref{fig:contractb} where the relevant parts are coloured in blue). This cobordism induces a homomorphism of vector spaces
\begin{eqnarray}
\label{multinH}
  \cF(W(c)b)\otimes \cF(W(b)a)\rightarrow\cF(W(c)a).
\end{eqnarray}

The algebra $\cH^n$ introduced in \cite{Khovanov} is defined as follows:
As a vector space it is
\begin{eqnarray}
  \cH^n=\bigoplus_{a,b\in\op{Cup}(n)}{}^{}_b\cH^n_a=\bigoplus_{a,b\in\op{Cup}(n)}\cF(W(b)a).
\end{eqnarray}
The elements from $\op{Cup}(n)$ should be thought of as being primitive idempotents of $\cH^n$, and the spaces $\cF(W(b)a)=:{}^{}_b\cH^n_a$ are the morphisms from the indecomposable projective left $\cH^n$-module indexed by $a$ to the one indexed by $b$. Therefore, one defines the product $fg=0$, if $f\in {}^{}_c\cH^n_d$, $g\in {}^{}_b\cH^n_a$, where $a,b,c,d\in\op{Cup}(n)$, $b\not=d$. In case $b=d$, the product is given by \eqref{multinH}.

If we consider $R$ as a graded vector space with the basis vector $1\in R$ in degree $-1$ and the basis vector $X\in R$ in degree $1$, then the vector space $\cH^n$ inherits a natural $\mZ$-grading. To make it compatible with the algebra structure we have to apply on overall shift $\langle n\rangle$ which increases the grading by $n$.  The graded vector space
$\cH^n=\bigoplus_{a,b\in\op{Cup}(n)}{}^{}_b\cH^n_a\langle n\rangle$
with the above multiplication becomes a positively graded algebra (\cite{Khotangles}). In the following we will consider the algebra $\cH^n$ with this grading. In particular, there are the subalgebras ${}^{}_a\cH^n_a\cong R^{\otimes n}\langle n\rangle$ of $\cH^n$ for any $a\in\op{Cup}(n)$ (cf. Lemma~\ref{Waa} below).

\subsection{Combinatorics: tableaux and generalised cup diagrams}
\label{combinatorics}

In this section we recollect a few combinatorial facts which will be needed later on. For any positive integer $n$ let $\cS(n)$ be the set of all sequences $\sigma=(\sigma_1, \sigma_2, \ldots, \sigma_{2n})$ where $\sigma_i\in\{+,-\}$ for $1\leq i\leq 2n$ with exactly $n$ pluses (minuses resp.). Of course, $S_{2n}$ acts transitively on $\cS(n)$ from the right hand side. Let $\sigma_{dom}:=(+,+,\ldots, +, -,-, \ldots, -)\in \cS(n)$.

Let $Y(n)$ be the set of Young diagrams which fit into an $n\times n$-square of $n^2$ boxes, i.e. the Young diagram corresponding to the partition $(n^n)$.  Let $Y(n, \op{upper})$ be the set of Young diagrams which fit into the Young diagram corresponding to the partition $(n-1, n-2, n-3,\ldots)$.

To $D\in Y(n)$ we associate a sequence $\sigma_D\in \cS(n)$ indicating its
shape as follows: first embed $D$ into a $n\times n$-square $D'$ of $n^2$
boxes, such that their upper left corners coincide. Now there is a unique path
$p_D$ from the bottom left corner to the top right corner of $D'$, along the
sides of the boxes, such that the interior of $D$ is completely to the left of
the path and all other boxes are to the right. The number of sides involved in
the path is always $2n$. Starting from the lower left corner, the path is
uniquely determined by giving the direction for each side. We use the rule
`minus=go up', `plus=go right'. In this way we associate to $D$ first a path $p_D$,
and then a sequence $\sigma_D\in\cS(n)$ encoding the path $p_D$.

\begin{ex}{\rm
It is $Y(2)=\left\{ {\tiny\Yv{\yng(2,2)}}\; , {\tiny\Yv\yng(2,1)}\; ,
{\tiny\Yv\yng(1,1)}\; , {\tiny\Yv\yng(2)}\; , {\tiny\yng(1)}\; ,
\emptyset\right\}$ and $Y(2,\op{upper})=\left\{{\tiny\yng(1)}\;
,\emptyset\right\}$. $\sigma_{\tiny\Yv{\yng(2,2)}}=(+,+,-,-)$,
$\sigma_{\tiny\Yv{\yng(2,1)}}=(+,-,+,-)$,
$\sigma_{\tiny\Yv{\yng(1,1)}}=(+,-,-,+)$,
$\sigma_{\tiny\Yv{\yng(2)}}=(-,+,+,-)$, $\sigma_{\tiny\Yv{\yng(1)}}=(-,+,-,+)$,
$\sigma_{\emptyset}=(-,-,+,+)$.}
\end{ex}

Let $\op{PMS}(n)$  (and $\op{PrInj}(n)$ respectively) be the set of iso-classes
of indecomposable projective(-injectives) modules in $\cO_0^{n,n}$.
\begin{prop}
\label{bijall}
  There are canonical bijections
  \begin{eqnarray*}
    \begin{array}[ht]{ccccccc}
      Y(n)&\leftrightarrow&\cS(n)&\leftrightarrow&S_n\times S_n\backslash S_{2n}&\leftrightarrow&\op{PMS}(n)\\
D&\mapsto&\sigma_D,\quad \sigma_{dom}w&\mapsto& w&\mapsto& [P^\p(w\cdot0)]
    \end{array}
  \end{eqnarray*}
\end{prop}

\begin{proof}
The first bijection is clear. The second is obvious, since $S_{2n}$ acts transitively on $\cS(n)$ and $\sigma_{dom}$ has stabiliser $S_n\times S_n$.  The third bijection is by definition (see Section~\ref{ordinary}).
  \end{proof}

To make the assignment $D\mapsto\sigma_D$ more precise we will now follow the
setup of Braden, see \cite{Braden} for details. Put $\mH:=\mZ+\frac{1}{2}$ and
for $k, l\in\mR$ set $\mH[k,l]:=\{\alpha\in \mH\mid k \leq\alpha\leq l\}$. We
generalise the construction above: Let $\la$ be a {\it partition}, by which
from now on we mean an infinite decreasing sequence
$\la=\la_1\geq\la_2\geq\la_3\cdots$  of non-negative integers such that
$\la_i=0$ for large $i$. We associate the corresponding Young diagram $D(\la)$
with $\la_i$ boxes in the $i$-row, and also an infinite $\{+,-\}$-sequence
$\varphi_\la$ indexed by $\mH$ (that is a function $\varphi_\la:
\mH\rightarrow\{-,+\}$) as follows:
\begin{eqnarray*}
\varphi_\la(\alpha)=
\begin{cases}
- &\text{if $\alpha=\la_i-i+\frac{1}{2}$ for some $i\in\mZ_{>0}$},\\
+&\text{otherwise.}
\end{cases}
\end{eqnarray*}
The uncommon indexing set is chosen to make it compatible with \cite{Braden}.
In particular, $D(\la)$ fits into a square of $n\times n$ boxes if and only if
\begin{itemize}
\item $\varphi_\la(j)=-$ if $j<-n+\frac{1}{2}$,
\item $\varphi_\la(j)=+$ if $j>n-\frac{1}{2}$,
\item The set $\{\varphi_\la(j)\mid n-\frac{1}{2}\geq j\geq -n+\frac{1}{2}\}$ contains exactly $n$ pluses and $n$ minuses.
\end{itemize}
Let $\widetilde{\cS(n)}$ be the set of such $\{+,-\}$-sequences. We have
isomorphisms of finite sets $Y(n)\cong \widetilde{\cS(n)}$,
$D(\la)\mapsto\varphi_\la$ and $\widetilde{\cS(n)}\cong \cS(n)$
$\varphi\mapsto(\varphi(-n+\frac{1}{2}), \varphi(-n+\frac{3}{2}), \ldots,
\varphi(n-\frac{1}{2}))$, the {\it restriction} of $\varphi$ to $\mH[-n, n]$.

\begin{definition}
\label{parent}
{\rm
Let $\la$ be a partition, i.e. $\la=\la_1\geq\la_2\geq\la_3\cdots$  of non-negative integers such that $\la_i=0$ for large $i$. Following \cite{Braden} we call
$(\alpha,\beta)\in\varphi_\la^{-1}(\{-\})\times\varphi_\la^{-1}(\{+\})$ {\it a
  $\la$-pair} if $\alpha<\beta$, $\sum_{\alpha\leq\gamma\leq\beta}\varphi_\la(\alpha)1=0$ and $\beta$ is minimal
with this property. If $(\alpha,\beta)$ and $(\alpha',\beta')$ are $\la$-pairs
then $(\alpha',\beta')>(\alpha,\beta)$ if and only if $\alpha'<\alpha$ and
$\beta'>\beta$. If $(\alpha',\beta')$ is minimal with this property  then
$(\alpha',\beta')$ is called a {\it parent} of $(\alpha,\beta)$. In the following $(\alpha',\beta')$ will always denote the parent of $(\alpha,\beta)$.
}
\end{definition}

For an example see the first diagram in Figure~\ref{fig:lanu}. Each $\la$-pair
has a unique parent (\cite[1.2]{Braden}), so the notation $(\alpha',\beta')$ makes sense. The set $\mH[-n, n]$ labels in a
natural way the endpoints of the arcs in any cup diagram from $\op{Cup}(n)$
from the left to the right: if $n=2$ for instance, then
$\operatorname{Cup}(2)=\{D,D'\}$ (see Figure~\ref{fig:cups}), the vertices
labelled by $-\frac{3}{2}$, $-\frac{1}{2}$, $\frac{1}{2}$, $\frac{3}{2}$ from
left to right. Then $(-\frac{3}{2},\frac{3}{2})$, $(-\frac{1}{2},\frac{1}{2})$
are $\la$-pairs for $D$, whereas $(-\frac{3}{2},-\frac{1}{2})$,
$(\frac{1}{2},\frac{3}{2})$ are $\la$-pairs for $D'$.\\

If $D\in Y(m,\op{upper})$ then let $c_D$ be the $2m$-cup diagram
$c_D\in\op{Cup}(m)$ where $\alpha, \beta\in\mH[-m, m]$ are connected if and
only if $(\alpha, \beta)$ form a $\la$-pair in $\varphi_\la$.  (This procedure
is well-defined since the pairs are nested (\cite[Lemma 1.2.1]{Braden}).

The bijections from the Proposition~\ref{bijall} provide three different
labelling sets for the isomorphism classes of indecomposable projective modules
in $\op{PMS}(n)$. The subset $\op{PrInj}(n)$ given by projective-injective modules
plays an important role, so we would like to have the corresponding labelling
sets singled out. The indecomposable projective module $P([w^\p]\cdot0)$
corresponding to the longest element $[w^\p]$ is always projective-injective in
$\cO^\p_0$, and all the other projective-injective modules are exactly the ones
which occur as direct summands when translation functors are applied to
$P([w^\p\cdot0])$, and their number equals the dimension of the Specht module
corresponding to $\la(\p)'$ (\cite{KMS}, Proposition~\ref{categorif}). Now
Proposition~\ref{bijall} restricts to
\begin{prop}
\label{bij}
  There are canonical bijections
  \begin{eqnarray*}
    \begin{array}[ht]{ccccccc}
      \op{Cup}(n)&\longleftarrow &Y(n,\op{upper})&\longleftrightarrow&\cS(n)'&\longleftrightarrow&\op{PrInj}(n)\\
c_D&\longmapsfrom&D&\longmapsto&\sigma_D,\quad\sigma_{dom}w&\longmapsto&
[P^\p(w\cdot0)]
    \end{array}
  \end{eqnarray*}
  where $\cS(n)'$ denotes the subset of sequences $\sigma\in\cS(n)$ such that for any plus there are
  more minuses than pluses appearing prior to
the given plus.
\end{prop}

  \begin{proof}
A sequence $\sigma\in\cS(n)$ can give rise to a cup diagram in $\op{Cup}(n)$ if
and only if for any plus there are more minuses than pluses appearing prior to
the given plus. Hence, the associated path stays above the diagonal, and gives
rise to a Young diagram $D\in Y(n,\op{upper})$ such that $\sigma_D=\sigma$.
Since the sets $\op{Cup}(n)$, $Y(n,\op{upper})$ and $\cS(n)'$, have the same
cardinality (the $n$-th Catalan number), the first two bijections follow. The
element $[w^\p]$ corresponds to the sequence $(-,-,\ldots, -,+,+,\ldots ,+)$,
hence to a $2n$-cup diagram. The last bijection follows from
Proposition~\ref{categorif}, since the $S_{2n}$-Specht module $S^{(n,n)}$
factors through the specialised Temperley-Lieb algebra $\op{TL}_{2n}$
(\cite[Theorem 4.1]{StDuke}) and the resulting representation is the regular
representation.
\end{proof}

\begin{ex} {\rm Let $n=1$, hence $\mg=\mathfrak{gl}_2$ with Weyl group $S_2=\langle s\rangle$ and $\p=\mb$.
The only element in $Y(1,\op{upper})$ is the empty diagram and corresponds to
the cup diagram with one cup and the sequence $(-,+)$, which then corresponds
to the indecomposable projective-injective module
$P^\p(s\cdot0)=P^\mb(s\cdot0)\in\cO_0^{1,1}$.}
\end{ex}

\begin{ex}{\rm
Let $n=2$, hence $\mg=\mathfrak{gl}_4$ with Weyl group $S_4=\langle
s_1=(1,2),s_2=(2,3),s_3=(3,4)\rangle$ and $\p$ such that $W_\p=\langle
s_1,s_3\rangle\cong S_2\times S_2$. The empty Young diagram corresponds to the
element $D$ from Figure~\ref{fig:cups} and the sequence $(-,-,+,+)$, which
correspond then to $P^\p(s_2s_1s_3s_2\cdot0)$, whereas the one box diagram
corresponds to $D'$ and the sequence $(-,+,-,+)$ which correspond to
$P:=P^\p(s_2s_1s_3\cdot0)$. Note that $P\cong\theta_\la^0P^\p(s_2\cdot\nu)$,
where $\nu$ is a dominant integral weight with stabiliser $\langle
s_1,s_3\rangle$ (hence illustrates Proposition~\ref{Irving}). It is easy to
check that $\END_\mg(P)\cong\mC[X]/(X^2)\otimes\mC[X]/(X^2)$.}
\end{ex}

\subsection{Khovanov's algebra and PrInj}
\label{KhovPrInj}

The motivation for the remaining sections is to establish a direct connection between the tangle invariants defined in \cite{Khotangles} on the one hand side and the ones defined in \cite{StDuke} on the other side.
In this paper we establish the key step from which it can then be deduced that Khovanov's invariants are nothing else than certain restrictions of the functorial invariants from \cite{StDuke} (see Section~\ref{outlook}).
The key step is to prove the following result (conjectured in \cite{StTQFT}, see also the weaker version in \cite{Braden}):
\begin{theorem}
\label{theresult}
  Let $n\in\mZ_{>0}$ and $\mg=\mathfrak{gl}_{2n}$. Let $\p=\p_n$ and $P^\p(x\cdot0)$, $x\in I\subseteq W^\p$ be a complete set of representatives for $\op{PrInj}(n)$. Set $$D_{n,n}:=\END_\mg(\oplus_{x\in I} P^\p(x\cdot0)).$$ Then there is an isomorphism of algebras
    \begin{eqnarray}
      D_{n,n}&\cong&\cH^n
    \end{eqnarray}
such that $\HOM_\mg(P^\p(x\cdot0),P^\p(y\cdot0))$ is identified with
${}^{}_b\cH^n_a$, where $a={\sigma_{dom}x}$, $b={\sigma_{dom}y}$. The isomorphism is even an isomorphism of $\mZ$-graded algebras.
\end{theorem}

\begin{cor}
  In the situation of the theorem we have $$\END_\mg(P)\cong(\mC[X]/(X^2))^{\otimes n}$$ for any indecomposable projective-injective module in $\cO_0^\p$.
\end{cor}

To prove the theorem we will embed the algebra $\cH^n$ into a larger algebra $\cK^n$ where the primitive idempotents are in bijection to the elements of $Y(n)$ and not just to the elements of $Y(n,\op{upper})$. The actual proof will be given in Section~\ref{SectionBraden}.

\subsection{The algebra  ${\cK^n}$, an enlargement of $\cH^n$}
\label{enlarged}

Let $a\in \cS(n)$. Take the corresponding Young diagram $D\in Y(n)$ (i.e. $a=\sigma_D$) and the corresponding partition $\la$. We view $D$ as a Young diagram $\tilde{D}\in Y(2n)$ and associate the $\{+,-\}$-sequence $\tilde{\sigma}_\la\in S(2n)$ of length $4n$ by restricting the $\{+,-\}$-sequence $\varphi_\la$ to $\mH[-2n, 2n]$. Alternatively, we could take the sequence $a$ and put $n$ minuses in front and $n$ pluses afterwards to obtain a $\{+,-\}$-sequence of length $4n$ which is exactly $\tilde{\sigma}_\la$.

Let $\tilde{a}=\sigma_{\tilde{D}}\in\op{Cup}(2n)$ be the corresponding $4n$-cup diagram where the arcs correspond to $\la$-pairs.
Hence, given $a$, $b\in S(n)$ we have the cup diagrams $\tilde{a}, \tilde{b}$ where the endpoints of the cups are labelled by $\alpha\in \mH[-2n,2n]$ and  the arcs correspond to $\la$-pairs. We call an endpoint {\it inner} if it is contained in  $\mH[-n, n]$, {\it outer left} if it is contained in  $\mH[-2n, -n]$, {\it outer right} if it is contained in $\mH[n, 2n]$.

Consider $W(\tilde{b})\tilde{a}$. This is a collection of circles which we view as coloured:
\begin{itemize}
\item A circle is {\it black} if it passes through inner points only.
\item A circle is {\it green} if it is not black and passes through at most one outer left point and at most one outer right point.
\item A circle is {\it red} if it is neither black nor green.
\end{itemize}

\begin{ex}{\rm For $n=2$ and $a=(+,+,-,-)$, $b=(-,+,+,-)$, $c=(+,-,+,-)$, $d=(-,-,+,+)$
we display in Figure~\ref{fig:colours} the cup diagrams $\tilde{a}$,
$\tilde{b}$, $\tilde{c}$ as well as the diagrams $W(\tilde{b})\tilde{a}$,
$W(\tilde{c})\tilde{a}$ and $W(\tilde{d})\tilde{d}$. (The dotted circles are
green; the left most (dashed) circle is red. There are two black circles.)}

\begin{figure}[htb]
  \centering
  \includegraphics{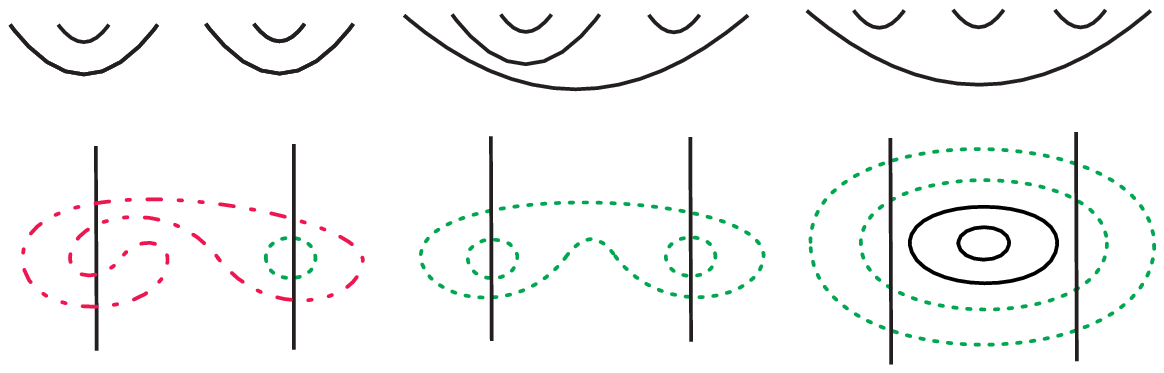}
  \caption{the vertical lines separate inner from outer points}
  \label{fig:colours}
\end{figure}

\end{ex}

Let ${\bf B}(b,a)$ (${\bf G}(b,a)$, ${\bf R}(b,a)$ respectively) be the number of black (green and red) circles in $W(\tilde{b})\tilde{a}$. To $W(\tilde{b})\tilde{a}$ we associate the complex vector space

$$\cG(W(\tilde{b})\tilde{a}):={}^{}_b\cK_a^n:=
\begin{cases}
  R^{\otimes{\bf B}(b,a)}\otimes \mC^{\otimes{\bf G}(b,a)}&\text{if ${\bf R}(b,a)=0$},\\
  \{0\}&\text{otherwise}.
\end{cases}
$$

As complex vector space, the algebra ${\cK^n}$ is
\begin{eqnarray}
  {\cK^n}=\displaystyle\bigoplus_{a,b\in\cS(n)}{}^{}_b\cK^n
_a=\bigoplus_{a,b\in\cS(n)}\cG(W(\tilde{b})\tilde{a}).
\end{eqnarray}
The unit
$\epsilon:\mC\rightarrow R$, the map $\epsilon':R\rightarrow\mC$ $1\mapsto 1$
$X\mapsto 0$, the inclusion $\{0\}\rightarrow R$ and the zero map
$R\rightarrow\{0\}$ give rise to canonical maps
$\op{can}:\cF(W(\tilde{c})\tilde{a}) \longrightarrow
\cG(W(\tilde{c})\tilde{a})$ and
$\op{can}:\cG(W(\tilde{c})\tilde{a})\rightarrow\cF(W(\tilde{c})\tilde{a})$
which `introduce the colouring' and `forget the colouring'. We turn ${\cK^n}$
into an algebra by putting $fg=0$ if $f\in {}^{}_c\cK^n_d$, $g\in
{}^{}_b\cK^n_a$, where $a,b,c,d\in \cS(n)$, $b\not=d$; and in case $b=d$, the
product is given by the composition
\begin{eqnarray}
\label{multinK}
  \begin{array}[ht]{ccc}
\cG(W(\tilde{c})\tilde{b})\otimes \cG(W(\tilde{b})\tilde{a})&&\cG(W(\tilde{c})\tilde{a}),\\
\downarrow&&\uparrow\\
\cF(W(\tilde{c})\tilde{b})\otimes \cF(W(\tilde{b})\tilde{a})&\rightarrow&\cF(W(\tilde{c})\tilde{a})
 \end{array}
\end{eqnarray}
where the vertical maps are the canonical ones and the horizontal one is the multiplication in $\cH^{2n}$ from \eqref{multinH}.

\begin{lemma}
\label{Waa}
  Let $a\in\cS(n)$. Then ${}^{}_a\cK^n_a\cong R^{\otimes{{\bf B}(a,a)}}\otimes\mC^{\otimes{{\bf G}(a,a)}}$ as algebras and the canonical map ${}^{}_{\tilde{a}}\cH^{2n}_{\tilde{a}}\rightarrow {}^{}_a\cK^n_a$ is surjective.
\end{lemma}

\begin{proof}
By definition of $\tilde{a}$ and the colouring rules, $W(\tilde{a})\tilde{a}$
is a union of black and green circles only. Hence ${}^{}_a\cK_a\cong
\cG(W(\tilde{a})\tilde{a})=R^{{\bf B}(a,a)}\otimes \mC^{\otimes{{\bf
G}(a,a)}}\not=\{0\}$. If we number the circles of $W(\tilde{a})\tilde{a}$,
first the black and then the green ones, then
\begin{eqnarray*}
\cG(W(\tilde{a})\tilde{a})\otimes \cG(W(\tilde{a})\tilde{a})&\longrightarrow&\cG(W(\tilde{a})\tilde{a}),\\
(R^{\otimes{{\bf B}(a,a)}}\otimes \mC^{\otimes{{\bf G}(a,a)}})\otimes (R^{\otimes{{\bf B}(a,a)}}\otimes \mC^{\otimes{{\bf G}(a,a)}})
&\longrightarrow&R^{\otimes{\bf B}(a,a)}\otimes \mC^{\otimes{{\bf G}(a,a)}}
\end{eqnarray*}
by multiplying the $i$-th factor in the first tensor product $\cG(W(\tilde{a})\tilde{a})$ with
the $i$-th factor in the second $\cG(W(\tilde{a})\tilde{a})$. The first statement of the lemma follows. The second follows directly from the definitions.
\end{proof}

\begin{ex}
\label{exampleSL2}
{\rm
  Let $n=1$ with the two sequences $a=(-,+)$ and $b=(+,-)$. Then $\tilde{a}$ and $\tilde{b}$ are the diagrams as depicted in Figure~\ref{fig:cups}. We get ${}^{}_a\cK^1_a=R\otimes\mC$ as algebra, ${}^{}_b\cK^1_a=\mC={}^{}_a\cK^1_b$ as vector spaces and ${}^{}_b\cK^1_b=\mC\otimes\mC$ as algebra. Using formula~\eqref{multinK} one easily verifies that the multiplication is given by the following formulas
  %in Table~\ref{multtable}
  %\begin{table}
  %\label{multtable}
  \begin{eqnarray*}
{}^{}_a\cK^1_a\otimes {}^{}_a\cK^1_b=R\otimes\mC\otimes\mC&\longrightarrow&\mC={}^{}_a\cK^1_b\\
1\otimes 1\otimes 1&\longmapsto& 1\\
X\otimes 1\otimes 1&\longmapsto&0.\\
 {}^{}_b\cK^1_a\otimes {}^{}_a\cK^1_b=\mC\otimes\mC&\longrightarrow&\mC\otimes\mC={}^{}_b\cK^1_b\\
1\otimes 1&\longmapsto& 0.\\
 {}^{}_a\cK^1_b\otimes {}^{}_b\cK^1_a=\mC\otimes\mC&\longrightarrow&R\otimes\mC={}^{}_a\cK^1_a\\
1\otimes 1&\longmapsto& X\otimes 1.\\
 {}^{}_a\cK^1_b\otimes {}^{}_b\cK^1_b=\mC\otimes\mC\otimes\mC&\longrightarrow&\mC={}^{}_a\cK^1_b\\
1\otimes 1\otimes 1&\longmapsto& 1.\\
 {}^{}_b\cK^1_a\otimes {}^{}_a\cK^1_a=\mC\otimes R\otimes\mC&\longrightarrow&\mC={}^{}_b\cK^1_a\\
1\otimes 1\otimes 1&\longmapsto& 1,\\
1\otimes X\otimes 1&\longmapsto&0.\\
 {}^{}_b\cK^1_b\otimes {}^{}_b\cK^1_a=\mC\otimes \mC\otimes\mC&\longrightarrow&\mC={}^{}_b\cK^1_a\\
1\otimes 1\otimes 1&\longmapsto& 1.
  \end{eqnarray*}
  %\caption{Multiplication in $\cK^1$}
  %\end{table}
Note that $\cK^1$ is isomorphic to the endomorphism ring $A_{1,1}$ of a minimal
projective generator in $\cO_0^\mb(\mathfrak{gl}_2)$ (\cite[5.1.1]{Stquiv} see
also Theorem~\ref{graphical}). }
\end{ex}

The $\mZ$-grading on $H^{2n}$ (see the paragraph before Section~\ref{combinatorics}) induces a unique grading on $K^n$ with respect to which both the maps $\CAN$ are graded maps. It is then clear that $K^n$ is itself a graded algebra. We leave it as an exercise to the reader to show that $K^1$ from Example~\ref{exampleSL2} is in fact a Koszul algebra.

If $a\in \cS(n)$ then we define for $1\leq i\leq 2n$ the element
$\op{can}(X_i(a))\in  {}^{}_a\cK^n_a$ as the image of $X_i(a):=1^{\otimes
  (i-1)}\otimes X\otimes 1^{\otimes (2n-i)}\in {}^{}_{\tilde{a}}\cH^{2n}_{\tilde{a}}$ under the canonical map.
The following lemma describes ${}^{}_b\cK^n_a$ as a left ${}^{}_b\cK^n_b$- and right ${}^{}_a\cK^n_a$-module.

\begin{lemma}
\label{rel4comb}
Let $a, b\in\cS(n)$. Let $f\in {}^{}_b\cK^n_a$ and $\op{can}(f)$ its canonical image in ${}^{}_{\tilde{b}}\cH^{2n}_{\tilde{a}}$. Then
\begin{eqnarray*}
  \op{can}(X_i(b))f&=&\op{can}\big(X_i(b)\op{can}(f)\big)\in {}^{}_b\cK^n_a\\
  f\op{can}(X_i(a))&=&\op{can}\big(\op{can}(f)X_i(a)\big)\in {}^{}_b\cK^n_a
\end{eqnarray*}
for $1\leq i\leq 2n$.
\end{lemma}

\begin{proof}
  This follows directly from the definitions \eqref{multinH} and \eqref{multinK}.
\end{proof}

\begin{figure}[htb]
  \centering
  \includegraphics[scale=1.1]{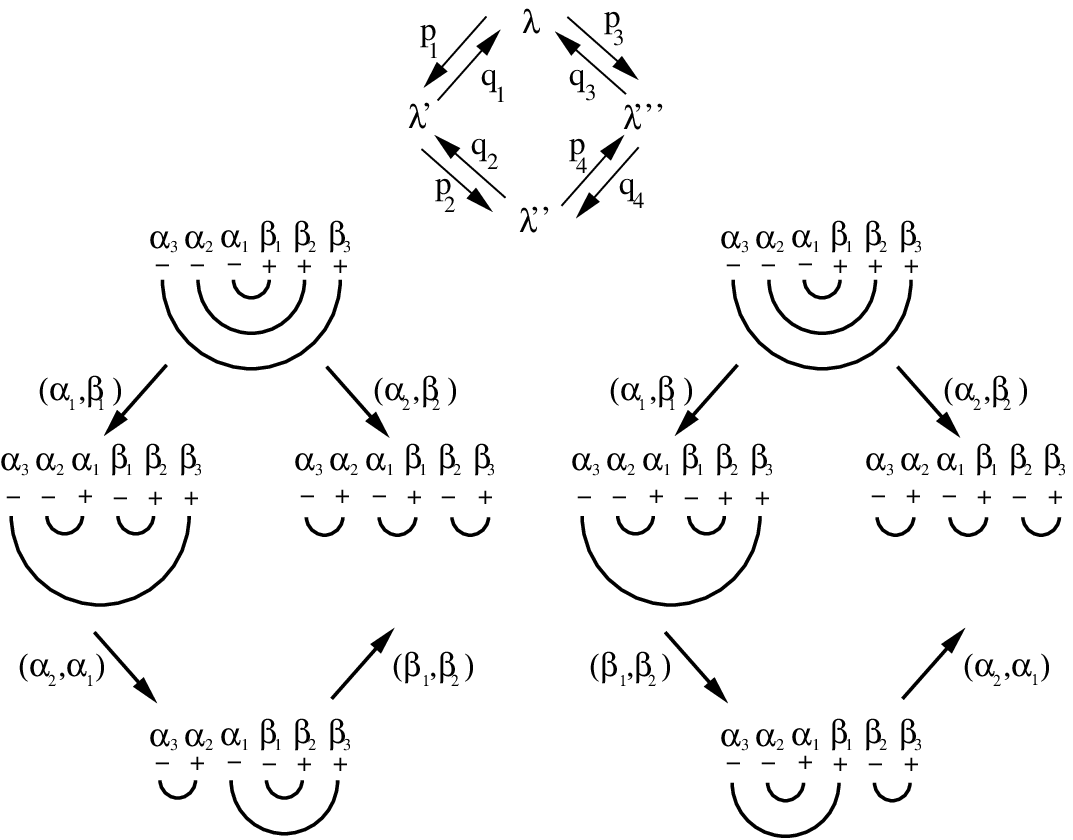}
  \caption{Diamonds}
  \label{fig:diamonds}
\end{figure}

\subsection{Braden's description of $\cO_0^{n,n}$}
\label{perverse}
If $P$ is a minimal projective generator of $\cO_0^{n,n}$ and $A_{n,n}=\END_\mg(P)$ then $\cO_0^{n,n}$ is equivalent to the category of finitely generated right $A_{n,n}$-modules. (In fact we could also work with left modules, since $\cO$ has a contravariant duality which identifies $A_{n,n}$ with its opposite algebra.) In \cite{Braden}, Braden gave an explicit description of $A_{n,n}$ in terms of generators and relations. We briefly recall this description.

\begin{definition}
{\rm
  Let $\la$ and $\nu$ be partitions. We write $\la\rightarrow\nu$ or $\la\stackrel{(\alpha,\beta)}\longrightarrow\nu$ if there is a $\la$-pair $(\alpha,\beta)$ such that $\varphi_\la(\gamma)=\varphi_\nu(\gamma)$ for any $\gamma\in\cH$, $\gamma\not=\alpha,\beta$ and $\varphi_\nu(\alpha)=+$ and $\varphi_\nu(\beta)=-$ (see Figure~\ref{fig:lanu}). Note that the $\la$-pair is uniquely determined by $\la$ and $\nu$. We write $\la\leftrightarrow\nu$ if either $\la\rightarrow\nu$ or $\nu\rightarrow\la$. Similarly, if $\phi_\la$, $\phi_\nu\in\widetilde{\cS(n)}$ with corresponding partitions $\la$ and $\nu$ then we write  $\phi_\la\rightarrow \phi_\nu$ if $\la\rightarrow\nu$ and  $\phi_\la\leftrightarrow \phi_\nu$ if $\la\leftrightarrow\nu$. A {\it diamond} is a tuple $(\la,\la',\la'',\la''')$ of four distinct partitions satisfying $\la\leftrightarrow\la'\leftrightarrow\la''\leftrightarrow\la'''\leftrightarrow\la$.
}
\end{definition}
Typical diamonds are depicted in Figure~\ref{fig:diamonds} where we display the relevant parts of the $\{+,-\}$-sequences with their cup-diagrams.

\begin{prop}(\cite[Section 1.3]{Braden})
\label{braden}
The algebra $A_{n,n}$ is the unitary associative $\mC$-algebra with generators
\begin{eqnarray*}
  &\left\{e_\la,t_{\alpha,\la}\mid\alpha\in\mH,\la\in \widetilde{\cS(n)}\right\},&\\
  &\left\{p(\la,\nu),\mu(\la,\nu)\mid \la,\nu\in \widetilde{\cS(n)}, \la\leftrightarrow\nu\right\},&
\end{eqnarray*}
and relations
\begin{enumerate}[(i)]
\item \label{rel1}$\sum_{\la\in \widetilde{\cS(n)}}e_\la=1$.
\item \label{rel2}$e_\la e_\nu=0$ if $\la\not=\nu$ and $e_\la e_\la=e_\la$ for any $\la\in \widetilde{\cS(n)}$.
\item \label{rel3}$\mu(\la,\nu)=1+p(\la,\nu)p(\nu,\la)$ for any $\la\leftrightarrow\nu$
\item \label{rel4}$t_{\alpha,\la}p(\la,\nu)=p(\la,\nu)t_{\alpha,\nu}$ for any $\alpha\in\mH$, $\la\in \widetilde{\cS(n)}$.
\item \label{rel5}$t_{\alpha,\la}t_{\beta,\nu}=0$ if $\la\not=\nu$.
\item \label{rel4b} The $t$'s commute with each other.
\item \label{rel6}$t_{\alpha,\la}t_{\beta,\la}=e_\la$ if $(\alpha,\beta)$ is a $\la$-pair.
\item \label{rel7}$t_{\alpha,\la}=e_\la$ for any $\la\in \widetilde{\cS(n)}$ if $\alpha<-n$ or $\alpha>n$.
\item \label{rel8}If $\la\stackrel{(\alpha,\beta)}\longrightarrow\nu$ and the $\la$-pair $(\alpha',\beta')$ is the parent of $(\alpha,\beta)$ then
  \begin{eqnarray*}
  \mu(\nu,\la)^{\eta(\beta)}&=&t_{\alpha,\nu}t_{\beta',\nu}\\
  \mu(\la,\nu)^{\eta(\beta)}&=&t_{\alpha,\la}t_{\beta',\la}
  \end{eqnarray*}
where $\eta(\beta)=(-1)^{\beta+\frac{1}{2}}$.
\item \label{rel9}If $(\la,\la',\la'',\la''')$ is a diamond with all elements contained in $\widetilde{\cS(n)}$ then
  \begin{eqnarray*}
    p(\la'',\la')p(\la',\la)&=&p(\la'',\la''')p(\la''',\la).
  \end{eqnarray*}
If all the elements in the diamond except $\la'''$ are in $\widetilde{\cS(n)}$ then
 \begin{eqnarray}
\label{zero}
    p(\la,\la')p(\la',\la'')&=0=&p(\la'',\la')p(\la',\la).
  \end{eqnarray}
\end{enumerate}
\end{prop}

\begin{remark}
\label{Schubertlabel}
{\rm The labeling of the idempotents in the algebra $A_{n,n}$ is such that $e_\la \in A_{n,n}=\END_\mg(P)$ is the idempotent projecting $P$ onto its summand $P^p(w.0)$ where $\la$ corresponds to $w$ according to Proposition~\ref{bij}.
For instance, the empty partition corresponds to the zero-dimensional Schubert cell, whereas the largest possible partition corresponds to the biggest Schubert cell. Indeed, the idempotent $e_\la$ is naturally associated with the Schubert cell $X_\la$ in \cite[1.1]{Braden}, which in turn comes along with an intersection homology complex $\cI_w$ corresponding to $P^\p(w\cdot0)$ (\cite[Theorem 12.2.5, also Example 12.2.6]{HTT}). The dictionary between the partition $\la$ and the Weyl group element $w$ is given by the formulas \cite[Proposition 8 and 9]{Fulton} together with the duality \cite[page 149]{Fulton}.
}
\end{remark}

\subsection{The map from $A_{n,n}$ to $\cK^n$}

The following lemma allows us to pass between $\la$-pairs, arcs and circles and follows directly from the definitions:

\begin{lemma}
\label{circles}
  Let $a\in\cS(n)$ and let $\la$ be the corresponding partition. Then
  \begin{enumerate}
  \item there is a canonical bijection between \\
the $\la$-pairs $(\alpha,\beta)$ where $-2n\leq\alpha\leq 2n$ and the cups in $\tilde{a}$;
  \item this bijection induces a canonical bijection between \\
the $\la$-pairs $(\alpha,\beta)$ where $-2n\leq\alpha\leq 2n$ and the circles in $W(\tilde{a})\tilde{a}$;
  \item  to any $\la$-pair $(\alpha,\beta)$  where $-2n\leq\alpha\leq 2n$ and $b\in\cS(n)$, there is a unique circle in $W(\tilde{b})\tilde{a}$ (and in $W(\tilde{a})\tilde{b}$ respectively) which contains the arc in $\tilde{a}$ corresponding to $(\alpha,\beta)$.
  \end{enumerate}
\end{lemma}

Let $a, b\in\cS(n)$ and $\la, \nu\in\widetilde S(n)$ their extensions. We denote by $e_\la, e_\nu$ the corresponding idempotents in $\cK^n$. If Lemma~\ref{circles} associates with $(\alpha,\beta)$ the $k$-th circle in $W(\tilde{a})\tilde{b}$ (or $W(\tilde{b})\tilde{a}$) then denote
\begin{eqnarray*}
  X_{\alpha}(a,b)&:=&\op{can}(1^{\otimes (k-1)}\otimes X\otimes 1^{\otimes 2n-k})\in {}^{}_a\cK^n_b, \\
X_{\alpha}(b,a)&:=&\op{can}(1^{\otimes (k-1)}\otimes X\otimes 1^{\otimes 2n-k})\in {}^{}_b\cK^n_a.
\end{eqnarray*}
If the pair $(\alpha,\beta)$ is not associated with a circle of $W(\tilde{a})\tilde{b}$ (or $W(\tilde{b})\tilde{a}$ respectively) then set $X_{\alpha}(a,b)=0\in{}^{}_a\cK^n_b$ (and $X_{\alpha}(b,a)=0\in{}^{}_b\cK^n_a$ resp.).

\begin{prop}
\label{arbeit}
In the notation of Definition~\ref{parent}, there is a homomorphism of algebras
\begin{eqnarray*}
  \begin{array}[thb]{cclr}
    \cE:\quad A_{n,n}&\longrightarrow&\cK^n\\
  e_\la&\longmapsto&
 e_\la&\in{}^{}_{a}\cK^n_{a}\\
  t_{\alpha,\la}&\longmapsto& e_\la+\eta(\beta)X_{\alpha}
&\in{}^{}_{a}\cK^n_{a}\\
  t_{\beta,\la}&\longmapsto& e_\la-\eta(\beta)X_{\alpha}&\in{}^{}_{a}\cK^n_{a}\\
p(\nu,\la)&\longmapsto& 1\otimes1\otimes\ldots\otimes 1+\frac{1}{2}X_{\alpha}&\in{}^{}_{b}\cK^n_{a}\\
p(\la,\nu)&\longmapsto& 1\otimes1\otimes\ldots\otimes 1+\frac{1}{2}X_{\alpha'}&\in{}^{}_{a}\cK^n_{b}\\
\mu(\la,\nu)&\longmapsto& e_\la+X_{\alpha}+X_{\alpha'}+X_{\alpha}\star X_{\alpha'}&\in{}^{}_{a}\cK^n_{a}
  \end{array}
\end{eqnarray*}
 where $\la$, $\nu\in\widetilde{\cS(n)}$ with restrictions $a$, $b\in\cS(n)$ and $\la$-pairs $(\alpha,\beta)$ such that $\la\stackrel{(\alpha,\beta)}\longrightarrow\nu$ in the last three cases; and $x\star y$ denotes the component-wise product in ${}^{}_a\cK^n_b$ for any $x$, $y\in {}^{}_a\cK^n_b$.
\end{prop}

For the proof we need the following lemma which is illustrated in Figure~\ref{fig:lanu}:

\begin{figure}[tbh]
  \centering
    \includegraphics{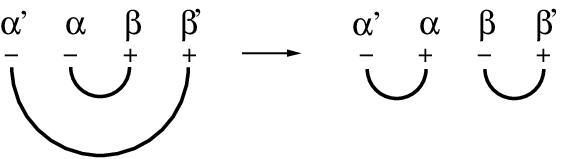}
  \caption{The relevant parts involved in $\la\stackrel{(\alpha,\beta)}\rightarrow\nu$}
  \label{fig:lanu}
\end{figure}

\begin{lemma}
\label{parentswap}
  Let $\la,\nu$ be partitions such that  $\la\stackrel{(\alpha,\beta)}\longrightarrow\nu$. Let $j':=(\alpha',\beta')$ be the parent of $j:=(\alpha,\beta)$. By definition of $\leftrightarrow$ the $\la$-pairs $j$ and $j'$ get transformed into the $\nu$-pairs $(\alpha',\alpha)$, $(\beta,\beta')$. Moreover, $\eta(\beta)=-\eta(\beta')$.
\end{lemma}

\begin{proof}
The definition of a parent (Definition~\ref{parent}) implies $\alpha'<\alpha$ and $\beta<\beta'$. The definition of $\leftrightarrow$ gives $\varphi_\nu(\alpha')=-$,  $\varphi_\nu(\alpha)=+$,  $\varphi_\nu(\beta)=-$,  $\varphi_\nu(\beta')=+$. Since $(\alpha',\beta')$ is the parent of $(\alpha,\beta)$, every element $x \in\mH[\alpha'+1,\alpha-1]$ must be $\la$-paired with an element $y\in\mH[\alpha'+1,\alpha-1]$. Therefore, $\sum_{\alpha'\leq\gamma\leq\alpha}\varphi_\nu(\gamma)1=\sum_{\alpha'<\gamma<\alpha}\varphi_\nu(\gamma)1=\sum_{\alpha'<\gamma<\alpha}\varphi_\la(\gamma)1=0$. Similarly  $\sum_{\beta\leq\gamma\leq\beta'}\varphi_\nu(\gamma)1=0$. Since going from $\la$ to $\nu$ only affects the pairs $j$ and $j'$ we are done. As all $x\in\mH[\beta,\beta']$ are $\nu$-paired inside this interval, the cardinality of $\mH[\beta,\beta']$ is even. In particular, $\eta(\beta)=-\eta(\beta')$.
\end{proof}

\begin{proof}[Proof of Proposition~\ref{arbeit}]
  We have to prove that the map $\cE$ is well-defined, hence to verify the compatibility with the relations from Proposition~\ref{braden}. Concerning the relations~\eqref{rel1}, \eqref{rel2} and \eqref{rel5}, there is nothing to do.

Let now $\la\stackrel{(\alpha,\beta)}\longrightarrow\nu$. Let
$j':=(\alpha'\beta')$ be the parent of $j:=(\alpha,\beta)$. Consider the
relation~\eqref{rel3}. Assume first that the circle corresponding to $j$ in
$W(\tilde{a})\tilde{b}$ as well as the circle corresponding to $j'$ in
$W(\tilde{b})\tilde{a}$ are both black and consider the corresponding subspaces
$R\subseteq{}^{}_{a}\cK^n_{b}$ and $R\subseteq{}^{}_{b}\cK^n_{a}$. Note that
for the description of the maps $\cE(p(\la,\nu))$ and $\cE(p(\nu,\la))$ only
the circles associated with $j$ and $j'$ are relevant. If we restrict ourselves
to the relevant circles, $\cE(p(\la,\nu))$ and $\cE(p(\nu,\la))$ become the
elements $1+\frac{1}{2}X\in R$ and their composition is displayed in the
Figures~\ref{fig:contract} and \ref{fig:contractb}, hence explicitly given as
\begin{eqnarray*}
\Delta\circ\op{m}
\big((1+\frac{1}{2}X) \otimes (1+\frac{1}{2}X)\big)
&=&\Delta\big((1+\frac{1}{2}X) \cdot (1+\frac{1}{2}X)\big)\\
=\quad\Delta(1+X)&=&X\otimes 1+1\otimes X+X\otimes X.
\end{eqnarray*}
Hence relation~\eqref{rel3} is satisfied if all the relevant circles are black. If there are green circles then the statement follows by applying the canonical maps. Amongst the relevant circles there are no red ones. To see this we consider again the diagrams in Figure~\ref{fig:contract} and Figure~\ref{fig:contractb}. If there is a red circle then at least one of the cup diagrams contains at least two left outer or two right outer points. This, however, is not possible because $a$, $b\in\cS(n)$.

 Therefore, $\cE$ is compatible with relation~\eqref{rel3}. The compatibility with relation~\eqref{rel4} is Lemma~\ref{rel4comb}. Relation~\eqref{rel4b} is clear.
If $(\alpha,\beta)$ is a $\la$-pair then
\begin{eqnarray*}
\cE(t_{\alpha,\la}t_{\beta,\la})
&=&\big(e_\la+\eta(\beta)X_{\alpha}\big)\big(e_\la-\eta(\beta)X_{\alpha}\big)  \\
&=&e_\la-\eta(\beta)X_{\alpha}+\eta(\beta)X_{\alpha}=e_\la=\cE(e_\la).
\end{eqnarray*}
Hence $\cE$ is compatible with relation~\eqref{rel6}. Relation~\eqref{rel7} holds by definition. Since $\mu(\la,\nu)$ is unipotent \cite[Proposition 1.8.2]{Braden} we could formally take the logarithm $\op{ln}$ of it and replace the relations~\eqref{rel8} by
\begin{eqnarray}
\label{log}
{\eta(\beta)}\op{ln}(\mu(\nu,\la))&=&\op{ln}(t_{\alpha,\nu})+\op{ln}(t_{\beta',\nu})\\
\label{log2}
{\eta(\beta)}\op{ln}(\mu(\la,\nu))&=&\op{ln}(t_{\alpha,\la})+\op{ln}(t_{\beta',\la}).  \end{eqnarray}
Note that for $Y=\eta(\beta)X_{\alpha}$ we have $\op{ln}(e_\nu+Y)=Y$. Therefore,
\begin{eqnarray*}
  \cE(\op{ln}(t_{\alpha,\nu})+\op{ln}(t_{\beta',\nu}))
&=&\eta(\beta)X_{\alpha}-\eta(\beta')X_{\alpha'}\\
&\stackrel{Lemma~\ref{parentswap}}=&\eta(\beta)X_{\alpha}+\eta(\beta)X_{\alpha'}.
\end{eqnarray*}
Now
$\cE(\op{ln}(\mu(\nu,\la)))=\op{ln}(e_\la+X_{\alpha}+X_{\alpha'}+X_{\alpha}\star
X_{\alpha'})=X_{\alpha}+X_{\alpha'}$, hence equation~\eqref{log},
and similarly \eqref{log2}, hold, and $\cE$ is compatible with the
relations~\eqref{rel8}. It is left to verify the diamond
relations~\eqref{rel9}. If $(\la,\la',\la'',\la''')$ is a diamond then the
$\varphi_\la$, $\varphi_{\la'}$, $\varphi_{\la''}$, $\varphi_{\la'''}$ agree
outside a set $N\subseteq\mH$ of cardinality $4$. On $N$, each of them takes
twice the value $-$ and twice the value $+$. Therefore, given $N$, there are
only six choices for $\la$. These, together with the possible diamonds, are
depicted in Figure~\ref{fig:diamonds}. Let us first assume that
\begin{eqnarray}
\label{out}
\quad\text{\it  All except the vertex $\la'''$ of the diagram are contained in $\widetilde{S(n)}$.}
\end{eqnarray}
In particular, there are either two left outer or two right outer points in $\la'''$ which are paired, and $\la'''$ is the only vertex in the diamond with this property. Let us consider the first diagram in Figure~\ref{fig:diamonds}. Let $A$, $B$, $C$, $D$ be the top, left, right, bottom vertex of the diamond respectively.

If $\la'''$ corresponds to $A$ then either $\alpha_1$ and $\beta_1$ are both left or both right outside, and then $\alpha_1, \alpha_2$ are left outside or $\beta_1,\beta_2$ are right outside in $B$. This contradicts~\eqref{out}.

If $\la'''$ corresponds to $B$ and say $\la$ to $A$ then $\la\stackrel{(\alpha_1,\beta_1)}\rightarrow\la'''$. Hence either $\alpha_1$ is outer and $\beta_1$ is inner or vice versa. In the first case $\alpha_3$ and $\alpha_2$ are both left outer in $D$, in the second case $\beta_1$ and $\beta_2$ are both right outer in $D$. This contradicts~\eqref{out}.

If $\la'''$ corresponds to $D$ and say $\la$ to $B$ then $\la\stackrel{(\alpha_2,\alpha_1)}\rightarrow\la'''$. Hence either $\alpha_2$ is outer and $\alpha_1$ is inner or vice versa. In the first case $\alpha_3$ and $\alpha_2$ are both left outer in $C$, in the second case $\beta_2$ and $\beta_3$ are both right outer in $C$. This contradicts~\eqref{out}.

If $\la'''$ corresponds to $C$ and say $\la$ to $A$ then $\la\stackrel{(\alpha_2,\beta_2)}\rightarrow\la'''$. Hence either $\alpha_2$ is outer and $\beta_2$ is inner or vice versa. In the first case $\alpha_3$ and $\alpha_2$ are both left outer in $D$. This contradicts~\eqref{out}. In the second case $W(A)D$ and $W(D)A$ are both a single circle which is red (since it connects $\beta_2$ and $\beta_3$). Hence
$$\cE(p(\la,\la'))\cE(p(\la',\la''))=0=\cE(p(\la'',\la'))\cE(p(\la',\la)).$$
The arguments for the second diamond depicted in Figure~\ref{fig:diamonds} are analogous. Hence $\cE$ is compatible with the relations~\eqref{zero}.

Let us now assume $(\la,\la',\la'',\la''')$ is a diamond with all elements contained in $\widetilde{\cS(n)}$. Let us consider again the diamonds in Figure~\ref{fig:diamonds}. For simplicity we first assume that all endpoints of the arcs are inner. Restricting to the relevant circles only we could consider any $\cE(p(\nu,\nu')$, $\nu\,\nu'\in\{\la,\la'\la''\la'''\}$ as an element of $R\otimes R$. An easy direct calculation shows that the composition of two of them is given by the multiplication
\begin{eqnarray}
\label{mapdiamonds}
  (R\otimes R)\otimes  (R\otimes R)&\rightarrow&R\\
   r_1\otimes r_2\otimes r_3\otimes r_4&\mapsto&r_1r_2r_3r_4
\end{eqnarray}
and the statement follows immediately.
If one of the relevant circles is red then not all elements of the diamond are contained in  $\widetilde{\cS(n)}$. If there are green circles appearing then the statement follows by applying the canonical map after and before ~\eqref{mapdiamonds}. (This is enough, because $\epsilon'\circ\epsilon: \mC\rightarrow R\rightarrow\mC$ is the identity map.) The map $\cE$ is therefore compatible with the relations~\eqref{rel9}.

So: the map $\cE$ is well-defined and gives rise to a homomorphism of algebras.
\end{proof}

\begin{prop}
\label{surjectivity}
  The algebra homomorphism $\cE$ is surjective.
\end{prop}

\begin{proof}
The algebra $\cK^n$ is by construction a graded quotient of the algebra $\cH^{2n}$.  Now, the algebra $\cH^{2n}$ is generated (over its semisimple degree zero part) in degrees one and two. To see this, first recall that for any $a\in\op{Cup}(2n)$, the subalgebra $a\cH^{2n}a$ is generated in degrees zero and two, whereas the space $a\cH^{2n}b$, equipped with its natural $a\cH^{2n}a$-module structure, is generated by its lowest degree element $1\in a\cH^{2n}b$. By \cite[Lemma1]{KhoSpringer}, the elements $1\in a\cH^{2n}b$ are contained in the subalgebra of $\cH^{2n}$ generated by degree one elements and so the claim follows. Now the proposition follows from Lemma~\ref{rel4comb} and the definition of the map $\cE$, since the image of $\cE$ contains the images of the generators of $\cH^{2n}$ in the quotient algebra $\cK^n$.
\end{proof}

\subsection{The grading of $A_{n,n}$ and the proof of Theorem~\ref{theresult}}
\label{SectionBraden}

Let $\la$, $\nu\in\widetilde{\cS(n)}$ and $\la\stackrel{(\alpha,\beta)}\longrightarrow\nu$ and consider $p:=p(\nu,\la)$ and $q:=p(\la,\nu)$. Set $x=qp$ and $y=pq$. One can find finite sums
$\tilde{p}:=\tilde{p}(\nu,\la):=\sum_{k\geq 0}c_kpx^k$ and $\tilde{q}:=\tilde{p}(\la,\nu):=\sum_{k\geq 0}c_kx^kq$, $c_k\in\mC$ such that $\tilde{p}\tilde{q}=\op{ln}(\mu(\nu,\la))$ and $\tilde{q}\tilde{p}=\op{ln}(\mu(\la,\nu))$. Namely define inductively $c_0=1$, and for $k\in\mZ_{>0}$ $$c_{k}=\frac{1}{2}\big((-1)^{k}\frac{1}{k+1}-\sum_{0<l,m<k,\;l+m=k}c_lc_m\big).$$
Then the claimed equalities hold (as formal power series in $x$ and $y$). Since $x$ and $y$ are nilpotent (\cite[Proposition 1.8.2]{Braden}) all the infinite sums are in fact finite. Note also that $1+x=\op{exp}(\op{ln}(1+x))$ is contained in the subalgebra generated by $1$ and $(\op{ln}(1+x))$, similarly for $y+1$. We get

\begin{prop}
\label{gradedBraden}
  The algebra $A_{n,n}$ is generated by
\begin{eqnarray*}
  &\left\{e_\la,\op{ln}(t_{\alpha,\la})\mid\alpha\in\mH,\la\in \widetilde{\cS(n)}\right\},&\\
  &\left\{\tilde{p}(\la,\nu),\op{ln}(\mu(\la,\nu))\mid \la,\nu\in \widetilde{\cS(n)}, \la\leftrightarrow\nu\right\}.&
\end{eqnarray*}
\end{prop}

\begin{proof}
  It is enough to show that the original generators are in the subalgebra, call it $B$, generated by the generators from the proposition. This is clear for the $e_\la$'s.
We also know it for the $\mu$'s. For the $t$'s it follows then from \cite[Proposition 1.8.1]{Braden}. Finally $\tilde{p}=p(1+\sum_{k>0}c_kx^k)$ where the second factor is invertible in $B$, since $x\in B$ is nilpotent. Hence $p\in B$ and similarly $q\in B$ and the statement follows.
\end{proof}

\begin{cor}
\label{corgrad}
  Putting the generators $e_\la$, $\op{ln}(t_{\alpha,\la})$, $\tilde{p}(\la,\nu),\op{ln}(\mu(\la,\nu))$ of $A_{n,n}$ from Proposition~\ref{gradedBraden} in degree $0$, $2$, $1$, $2$ turns $A_{n,n}$ into a positively graded algebra and $\cE$ becomes a homomorphism of $\mZ$-graded algebras. This grading is the Koszul grading.\end{cor}

\begin{proof}
Using the new generators, the relations from Proposition~\ref{braden} become
homogeneous. This is completely obvious except for relation~\eqref{rel9}. Let
us assume this to be true for the moment then the relations also show that
$A_{n,n}$ becomes a quadratic positively graded algebra, i.e. generated in
degree zero and one with relations in degree two. Now we are in the situation
of \cite[Proposition 2.4.1]{BGS}, i.e. for any graded $A_{n,n}$-module $M$ with
simple head, the radical filtration of $M$ agrees (at least up to a shift in
the grading) with the grading filtration. This holds in particular for
indecomposable projective modules. The same holds if we equip $A_{n,n}$ with
its Koszul grading (\cite{BGS}). By the unicity of gradings (\cite[2.5]{BGS}),
the statement follows if we verified relation~\eqref{rel9}.

Since we only have a case by case argument, we will sketch the argument for a specific example only and leave it to the reader to figure out all other possibilities.
Let $(\la,\la',\la'',\la''')$ be the diamond on the left in Figure~\ref{fig:diamonds}. We use the notation from the small diamond in the middle of Figure~\ref{fig:diamonds}. We claim that the relation $p_2p_1=q_4p_3$ could be replaced by the relation $\tilde{p}_2\tilde{p}_1=\tilde{q}_4\tilde{p}_3$. Let $x_i=q_ip_i$. Using the relations of Proposition~\ref{braden} we get
\begin{eqnarray*}
 p_2p_1(1+x_1)&\stackrel{\eqref{rel3}, \eqref{rel8}}=&p_2p_1(t_{\alpha_1,\la}t_{\beta_2,\la})^{\eta(\beta_1)}\stackrel{\eqref{rel4}}=(t_{\alpha_1,\la''}t_{\beta_2,\la''})^{\eta(\beta_1)}p_2p_1\nonumber\\
&\stackrel{\eqref{rel6}, \eqref{rel4b}}=&(t_{\beta_1,\la''}t_{\beta_3,\la''})^{-\eta(\beta_1)}p_2p_1\stackrel{\eqref{rel4}}=p_2p_1T\nonumber\\
p_2(1+x_2)p_1&\stackrel{\eqref{rel3}, \eqref{rel8}}=&p_2(t_{\alpha_2,\la'}t_{\beta_2,\la'})^{\eta(\alpha_1)}p_1\stackrel{\eqref{rel4}}=p_2p_1T',\nonumber\\
q_4p_3(1+x_3)&\stackrel{\eqref{rel3}, \eqref{rel8}}=&q_4p_3T'\nonumber\\
(1+x_4)q_4p_3&\stackrel{\eqref{rel3}, \eqref{rel8}}=&(t_{\beta_1,\la'''}t_{\beta_3,\la'''})^{\eta(\beta_2)}q_4p_3\stackrel{\eqref{rel4}}=q_4p_3T.
\end{eqnarray*}
where $T:=(t_{\beta_1,\la}t_{\beta_3,\la})^{-\eta(\beta_1)}$ and $T'=(t_{\alpha_2,\la}t_{\beta_2,\la})^{\eta(\beta_2)}$.

Hence $p_2p_1=q_4p_3$ implies $p_2(1+x_2)^k p_1=q_4p_3(1+x_3)^k$ for any $k\in\mZ_{\geq0}$ and then $p_2(x_2)^kp_1=q_4p_3(x_3)^k$ by induction. The displayed equations above also imply $p_2(1+x_2)^k p_1(1+x_1)^l=(1+x^4)^lq_4p_3(1+x_3)^k$ for any $k,l\in\mZ_{\geq0}$ and then by induction $p_2(x_2)^kp_1(x_1)^l=(x_4)^lq_4p_3(x_3)^k$. It follows $\tilde{p}_2\tilde{p}_1=\tilde{q}_4\tilde{p}_3$.

Assume now that $\tilde{p}_2\tilde{p}_1=\tilde{q}_4\tilde{p}_3$. Then $\tilde{p}_2(1+x_2)^k \tilde{p}_1=\tilde{q}_4\tilde{p}_3(1+x_3)^k$ for any $k\in\mZ_{\geq0}$ with the arguments from above. By induction we obtain  $\tilde{p}_2(x_2)^k \tilde{p}_1=\tilde{q}_4\tilde{p}_3(x_3)^k$. Using again the arguments from above we get $\tilde{p}_2(x_2)^k \tilde{p}_1(1+x_1)^l=(1+x_4)^l\tilde{q}_4\tilde{p}_3(x_3)^k$ for any $k,l\in\mZ_{\geq0}$. By induction we therefore get
\begin{eqnarray}
\label{tilde}
   \tilde{p}_2(x_2)^k\tilde{p}_1(x_1)^l=(x_4)^l\tilde{q}_4\tilde{p}_3(x_3)^k.
\end{eqnarray}
Now let $l$ be maximal such that $p_1x_1^l\not=0$ or $x_4^lq_4\not=0$. Then choose (if possible) $k$ maximal such that $p_2(x_2)^kp_1(x_1)^l\not=0$ or $(x_4)^lq_4p_3(x_3)^k\not=0$. From \eqref{tilde} we get
\begin{eqnarray}
  \label{eq:co}
  c_0c_0p_2(x_2)^kp_1(x_1)^l&=&c_0c_0(x_4)^lq_4p_3(x_3)^k.
\end{eqnarray}
Then we choose $k'<k$ maximal with the above conditions and deduce that $p_2x_2^{k'}p_1x_1^l=x_4^lq_4p_3x_3^{k'}$. Inductively, the latter holds for any $k'$. By double induction on $l$ and $k$ we finally obtain $p_2p_1=q_4p_3$.
\end{proof}

We get a description of the arrows in the Ext-quiver of the algebra $A_{m,m}$:

\begin{cor}
  \label{Ext}
Let $L(v\cdot0)$, $L(w\cdot0)\in\cO_0^{n,n}$ be simple modules. Consider $\sigma_{dom}v$, $\sigma_{dom}w\in\cS(n)$ and let $\nu$, $\la\in\widetilde{\cS(n)}$ be their  extensions. Then
  \begin{eqnarray*}
    \EXT^1_{\cO_0^{n,n}}(L(v\cdot0), L(w\cdot0))=
    \begin{cases}
    \mC&\text{if $\nu\leftrightarrow\la$}\\
    \{0\}&\text{otherwise}.
    \end{cases}
  \end{eqnarray*}
\end{cor}

\begin{proof}
  Let $P^{(n,n)}(v\cdot0)$, $P^{(n,n)}(w\cdot0)\in\cO_0^{(n,n)}$ be the projective cover of
  $L(v\cdot0)$ and $L(w\cdot0)$ respectively.
Then the dimension of $\EXT^1_{\cO_0^{n,n}}(L(v\cdot0),
  L(w\cdot0))$ is equal to the dimension of the subspace $M$ of
  $\HOM_{\cO}(P^{(n,n)}(w\cdot0),P^{(n,n)}(v\cdot0))$ spanned by all morphisms
  $f$ whose image is contained in the radical of $P^{(n,n)}(v\cdot0)$, but not
  in the square of the radical. When passing to $A_{n,n}$, this subspace
  corresponds to the space spanned be all morphisms of degree one, since, for indecomposable projectives, the grading filtration
  agrees with the radical filtration thanks to \cite[Proposition 2.4.1]{BGS}. By
Corollary~\ref{corgrad} there is, up to a nonzero scalar, a unique morphism of
degree one if  $\nu\leftrightarrow\lambda$ and no morphism otherwise.
\end{proof}

\begin{proof}[Proof of Theorem~\ref{theresult}]
The homomorphism $\cE$ from Proposition~\ref{arbeit} induces a surjective homomorphism of algebras
    \begin{eqnarray}
      \cE':\quad D_{n,n}&\rightarrow&\cH^n
    \end{eqnarray}
such that $\HOM_\mg(P^\p(x\cdot0),P^\p(y\cdot0))$ is mapped to $_a\cH^n_b$,
where $a={\sigma_{dom}x}$, $b={\sigma_{dom}y}$ (see Remark~\ref{Schubertlabel}).
By Corollary~\ref{corgrad} it
is only left to show that $\cE'$ is an isomorphism. To do so it is enough to
compare the dimensions. However, the dimension of $_{a}\cH^n_{b}$ is $2^k$,
where $k$ is the number of circles in $W(a)b$. We could rephrase this as
follows: Consider the irreducible right (complex) $S_{2n}$-module $M$
corresponding to the partition $2n=n+n$. This module has a unique up to a
scalar symmetric non-degenerate $S_{2n}$-invariant bilinear form $\bf{b}$
(\cite[Section 6]{Murphy}). One can naturally identify the elements of
$\op{Cup}(n)$ with the basis of $M$ obtained by specialising the
Kazhdan-Lusztig basis in the generic Hecke algebra of $S_{2n}$ such that
$\DIM_\mC(_a\cH^n_b)={\bf b}(a,b)$ (\cite[Theorem 7.3]{Fung} and references
therein.) On the other hand, we categorified $M$ in
Proposition~\ref{categorif}. One can also categorify the bilinear forms as
follows: there is a scalar $\gamma\in\mC$ such that
\begin{eqnarray}
\label{bilform}
\DIM_\mC\HOM_\mg(P^\p(x\cdot0),P^\p(y\cdot0))={\bf
b}(a,b)\gamma
\end{eqnarray}
where $a={\sigma_{dom}x}$, $b={\sigma_{dom}y})$ (\cite[Proposition 4]{KMS}).
The explicit formulas in \cite[Corollary p.327]{IS} give the dimension of the
endomorphism ring of $P\in\cO_0^{n,n}$ as in Proposition~\ref{Irving}, namely
as follows: $P\cong \theta_\nu^0 M^\p(\nu)$ where $M^\p(\nu)$ is a (simple,
projective) generalised Verma module. By \cite[Proposition 3.1]{IS} the
stabiliser of $\nu$ is isomorphic to $n$ copies of $S_2$. This implies
$\DIM_\mC\HOM_\mg(P,P)=2^n$ and hence $\gamma=1$. The theorem follows.
\end{proof}

\subsection{A graphical description of $\cO^{n,n}_0$}

In this section we prove the following

\begin{theorem}
\label{graphical}
  For any $n\in\mZ_{>0}$ the algebra homomorphism $\cE: A_{n,n}\cong\cK^n$ is an isomorphism.
\end{theorem}

\begin{proof}
We use the notation from Proposition~\ref{bijall} and Section~\ref{enlarged} and set $\p=\p_n$. For $w\in W^\p$ put $a_w=\sigma_{dom}w$. Let us denote by $\varphi_w\in\widetilde{\cS(n)}$ the extension of $a_w$, and let $\la_w$ be the corresponding partition.

Since $\cE$ is surjective (Proposition~\ref{surjectivity}) and the involved algebras are finite dimensional, it is enough to show that for any $v,w\in W^{\p}$
  \begin{eqnarray}
\label{claim}
    \HOM_\mg\big(P^{\p}(v\cdot0),P{^\p}(w\cdot0)\big)&\cong&\cG\big(W(\tilde{a}_v)\tilde{a}_w\big).
  \end{eqnarray}

Let us first consider the case where $\nu=e$, the identity of $W$. By \cite[Corollary 5.2]{Brenti} we have $\{0\}\not=\HOM_\mg\big(P^\p(0),P^\p(w\cdot0)\big)$ (and then equal to $\mC$) if and only if $${a}_w=
(\underbrace{+,\ldots, +}_{r},
\underbrace{-\ldots -}_{s},
\underbrace{+,\ldots +}_{s},
\underbrace{-\ldots -}_{r})=:+^r-^s+^s-^r$$ for some $r,s\in\mZ_{\geq0}$.

If there are $\gamma, \delta\in\mH[-n,0]$ $\gamma<\delta$ such that $\varphi_w(\gamma)=-$ and $\varphi_w(\delta)=+$ then there is a $\la_w$-pair $(\alpha, \beta)$ such that $\alpha$, $\beta\in\mH[-n,0]$. Moreover, $\varphi_e(\alpha)=\,+\,=\varphi_e(\beta)$. The corresponding circle in $W(\tilde{a}_w)\tilde{a}_e$ is then red. Similarly, if there are $\gamma, \delta\in\mH[0,n]$ $\gamma<\delta$ such that $\varphi_w(\gamma)=-$ and $\varphi_w(\delta)=+$ then there is a $\la_w$-pair $(\alpha, \beta)$ such that $\alpha$, $\beta\in\mH[0,n]$, $\varphi_e(\alpha)=\,-\,=\varphi_e(\beta)$. The corresponding circle in $W(\tilde{a}_w)\tilde{a}_e$ is red. Therefore $\cG\big(W(\tilde{a}_w)\tilde{a}_e\big)=\{0\}$ if $a_w$ is not of the form $+^r-^s+^s-^r$ as above. On the other hand if $a_w=+^r-^s+^s-^r$ then  $W(\tilde{a}_w)\tilde{a}_e$ consists of green circles only (as depicted in Figure~\ref{fig:Brentilemma}), and  $\cG\big(W(\tilde{a}_v)\tilde{a}_w\big)=\mC$.

\begin{figure}[bth]
  \centering
   \includegraphics[scale=0.5]{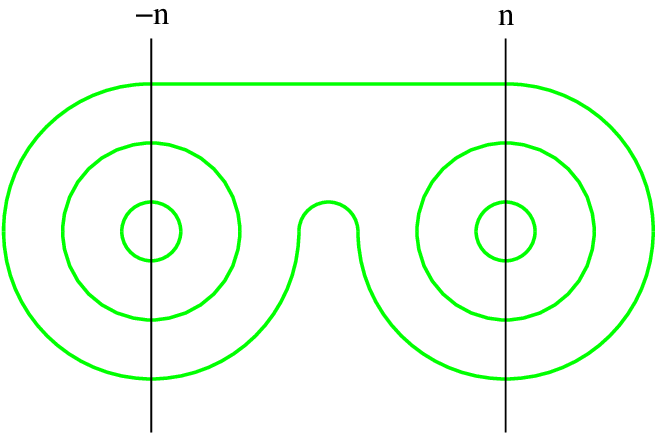}
  \caption{$W(\tilde{a}_w)\tilde{a}_e$ consists of green circles only.}
  \label{fig:Brentilemma}
\end{figure}

Hence formula~\eqref{claim} is true for $v=e$ and we can do induction on the length of $v$. Choose some simple reflection $s$ such that $vs<v$. Then
\begin{eqnarray*}
&&   \HOM_\mg\big(P^{\p}(v\cdot0),P{^\p}(w\cdot0)\big)\\
&=&\HOM_\mg\big(\theta_s P^{\p}(vs\cdot0),P{^\p}(w\cdot0)\big)\\
&=&\HOM_\mg\big(P^{\p}(vs\cdot0),\theta_sP{^\p}(w\cdot0)\big)\\
&=&
\begin{cases}
  \HOM_\mg\big(P^{\p}(vs\cdot0),P{^\p}(ws\cdot0)\big)&\text{$ws\in W^{\p}$, $ws>w$,}\\
  \HOM_\mg\big(P^{\p}(vs\cdot0),P{^\p}(w\cdot0)\oplus P{^\p}(w\cdot0)\big)&\text{if $ws<w$,}\\
\{0\}&\text{if $ws\not\in W^\p$}.
\end{cases}
\end{eqnarray*}
by \cite[Theorem 1]{BW}, the self-adjointness of $\theta_s$ and \cite[Section
3]{SoKipp}. On the other hand, if $ws\in W^{\p}$, $ws>w$ then obviously
$\cG\big(W(\tilde{a}_v)\tilde{a}_w\big)=\cG\big(W(\tilde{a}_{vs})\tilde{a}_{ws}\big)$
which is zero if $ws\notin W^{\p}$. If $ws\in W^{\p}$ and $ws<w$ then
$\cG\big(W(\tilde{a}_v)\tilde{a}_{ws}\big)$ differs from
$\cG\big(W(\tilde{a}_v)\tilde{a}_w\big)$ by a black circle. Formula
\eqref{claim} follows therefore from the induction hypothesis.
\end{proof}

\begin{remark}{\rm
  Theorem~\ref{graphical} can be
  generalised to all maximal parabolic
  subalgebras or perverse sheaves on Grassmannians $\op{Perv}_B(Gr(k,n))$.
  In general one has to take $\{+,-\}$-sequences of length $n$ with $k$ pluses and $n-k$ minuses. They define the inner points.
  Then we add $n-k$ minuses to the left and $n-(n-k)$ pluses to the right and proceed as before.
}
\end{remark}

\subsection{Intersections of components of the Springer fibre}
In the first three Sections of the paper we used Soergel's Endomorphismensatz
(\cite{Sperv}) which in particular implies that the endomorphism ring of the
only indecomposable projective-injective module in $\cO_0$ has commutative
endomorphism ring, isomorphic to $H^*(\cB)$ (Proposition~\ref{centreSoerg} and
the preceding paragraph). On the other hand we know by
Proposition~\ref{maxparab} that the endomorphism ring of any indecomposable
projective module in a $\cO^\p$ for maximal parabolic $\p$ is commutative, and
Theorem~\ref{theresult} gives an explicit algebraic description of the
endomorphism rings of indecomposable projective-injective modules. In fact it
describes the space of homomorphisms between two indecomposable
projective-injective modules as a bimodule over their endomorphism rings. A
geometrical interpretation of these bimodules in terms of cohomology rings is
still missing. We would like to finish this paper by formulating a conjectural
interpretation based on the following:

\begin{theorem}
\label{PrInjSpringer} Let $\mg=\mathfrak{gl}_n$ and let $\p$ be some maximal
parabolic subalgebra. Let $\op{Irr}(\cB_\p)$ denote the set of irreducible
components of $\cB_\p$ and  $\op{PrInj}(\p)$ the set of isomorphism classes of
indecomposable projective-injective modules in $\cO_0^\p$. Then there is a
bijection $$\psi:\quad\quad \op{PrInj}(\p)\cong\op{Irr}(\cB_\p)$$ such that
there is an isomorphism of vector spaces $$\HOM_\mg(P,Q)\cong
H^*(\psi(P)\cap\psi(Q))$$ for any $P$, $Q\in\op{PrInj}(\p)$.
\end{theorem}

\begin{proof}
By results of Vargas and of Spaltenstein, there is an explicit bijection
between the irreducible components of $\cB_\p$ and standard tableaux of shape
$\la(\p)$. In our special situation we have to consider only tableaux with two
rows, and we refer to \cite[Theorem 5.2]{Fung}, where the bijection is made
explicit. We restrict ourselves to the case where $\mu=\la(\p)$ is already a
partition. This is possible thanks to \cite[Theorem 5.4]{MSSerre}.

Let us first consider the case $\la(\p)=(n,n)$. Given a standard tablaux $T$
with two rows of length $\la_1=n$ and $\la_2=n$, one can associate a cup
diagram in $\operatorname{Cup}(n)$ as follows: $T$ has the entries
${1,2,\ldots, 2n}$ so that the numbers are decreasing from left to right in
each row, and decreasing from top to bottom in each column. The cup diagram
$C_T$ has $2n$ vertices, labelled by $1$ to $2n$ from the left, so that the
left endpoint of each cup is labelled by a number appearing in the bottom row
of $T$, whereas the right endpoints of the cups are labelled by elements from
the top row of $T$. (So the end points of each cup are in different rows of
$T$). For example, if $n=2$ then we have the standard tableaux
$$\young(4,2)\young(3,1)\quad\quad\young(4,3)\young(2,1)$$ to which we associate
$D$ and $D'$ displayed in Figure~\ref{fig:cups}. It is easy to see that this
procedure provides a bijection between standard tableaux of shape $(n,n)$ (and
hence of irreducible components of $\cB_\p$), and cup diagrams from
$\operatorname{Cup}(n)$. Now we apply Proposition~\ref{bij} and obtain a
bijection between the irreducible components and the isomorphism classes of
indecomposable projective-injective modules
 $\operatorname{PrInj}(n)$. Thanks to Theorem~\ref{theresult} the dimension of the homomorphism spaces
 between projective-injective modules can be computed using the diagram
 calculus, which is also used in \cite[Theorem 7.2, Theorem 7.3]{Fung} to compute the dimension of the
cohomology of the intersection of two components. This settles the case of the
partition $(n,n)$ and gives an explicit way to compute the dimensions of the
vector spaces.

In general, one should argue as follows: first of all we have
(\cite[Proposition 4]{KMS}) the categorification of the Specht modules for the
symmetric group (Proposition~\ref{categorif}), but also of its invariant
bilinear form \eqref{bilform} by taking dimensions of the homomorphism spaces
between indecomposable projective modules. On the other hand, the dimension of
the cohomology of the intersection of components is computed in the same way
\cite[Theorem 7.2, Theorem 7.3]{Fung}, so that the statement follows up to a
multiplication with a common factor. But one easily checks that this common
factor must be equal to one by computing one of the endomorphism rings
explicitly.
\end{proof}

Up to a shift, the isomorphism from the theorem is an isomorphism of
$\mZ$-graded vector spaces, where the grading on the hom-space is given as in
Section~\ref{gradings}.

We conjecture that this isomorphism is compatible with the multiplicative
structure as well:

\begin{conjecture}
\label{conj}
  The isomorphism $\END_\mg(P)\cong H^*(\psi(P))$ is a ring homomorphism and  $$\HOM_\mg(P,Q)\cong H^*(\psi(P)\cap\psi(Q))$$ as $\big(\END_\mg(P), \END_\mg(Q)\big)=(H^*(\psi(P),H^*(\psi(Q))$-module.
\end{conjecture}

\subsection{Connection with Khovanov homology}
\label{outlook}
Using Theorem~\ref{theresult} one can deduce that the full Conjecture \cite[Conjecture 2.9]{StTQFT} holds. Roughly speaking this says that the functorial tangle invariants defined in \cite{Khotangles} are obtained by restring the functorial tangle invariants from \cite{StDuke} to a certain subcategory, invariant under all these functors. In particular, the resulting homological invariants are the same. The proof is quite technical and lengthy and will appear in a subsequent paper. We expect that Conjecture~\ref{conj} provides the basis for a geometric interpretation of Khovanov homology using categories of sheaves related to the Springer fibres.

\bibliography{ref}
%Included for Gather Purpose only:
%input "ref.bib"
\end{document}